\documentclass[12pt]{amsart}
\usepackage{a4wide,enumerate,xcolor}
\usepackage{amsmath,graphicx,comment,amssymb}
\allowdisplaybreaks
\usepackage[utf8]{inputenc}
\usepackage[english,vietnamese]{babel}

\usepackage{enumitem}

\let\pa\partial
\let\na\nabla
\let\eps\varepsilon
\newcommand{\N}{{\mathbb N}}
\newcommand{\R}{{\mathbb R}}
\newcommand{\diver}{\operatorname{div}}

\newtheorem{theorem}{Theorem}[section]
\newtheorem{lemma}[theorem]{Lemma}
\newtheorem{proposition}[theorem]{Proposition}
\newtheorem{remark}[theorem]{Remark}

\begin{document}

\title[Three-species drift-diffusion model for memristors]{Time-uniform bounds and long-time asymptotics for a three-species drift-diffusion model for memristors}

\author[A. J\"ungel]{Ansgar J\"ungel}
\address{Institute of Analysis and Scientific Computing, TU Wien, Wiedner Hauptstra\ss e 8--10, 1040 Wien, Austria}
\email{juengel@tuwien.ac.at} 

\author[\selectlanguage{vietnamese}T. T. Nguyễn]{ \selectlanguage{vietnamese}Tuấn Tùng Nguyễn}
\address{Institute of Analysis and Scientific Computing, TU Wien, Wiedner Hauptstra\ss e 8--10, 1040 Wien, Austria}
\email{tung.nguyen@tuwien.ac.at}

\selectlanguage{english}
\date{\today}

\thanks{The authors acknowledge partial support from the Austrian Science Fund (FWF), grant 10.55776/PAT2687825, and from the Austrian Federal Ministry for Women, Science and Research and implemented by \"OAD, grant MULT09/2025. This work has received funding from the European Research Council (ERC) under the European Union's Horizon 2020 research and innovation programme, ERC Advanced Grant NEUROMORPH, no.~101018153. For open-access purposes, the authors have applied a CC BY public copyright license to any author-accepted manuscript version arising from this submission.} 

\begin{abstract}
The qualitative behavior of weak solutions to the two-dimensional drift-diffusion equations for memristors, incorporating Shockley--Read--Hall-type generation--recombination terms and time-dependent boundary conditions, is investigated. The model couples the electron, hole, and oxygen-vacancy densities to the electric potential through the Poisson equation. The global existence and uniqueness of weak solutions are established. If the boundary data converge to the thermal equilibrium at an integrable rate, the densities remain uniformly-in-time bounded above and below, and the solution converges exponentially fast to thermal equilibrium. For time-independent boundary data and sufficiently small applied bias, solutions converge exponentially fast to steady states carrying nonzero fluxes.
\end{abstract}

% \paragraph{Keywords:}  

% \paragraph{AMS classification:}  

\maketitle

\section{Introduction}

Neuromorphic computing is a brain-inspired computational paradigm that seeks to enhance the energy efficiency of artificial neural networks through dedicated hardware architectures and algorithms. Its importance is growing in response to the increasing computational demands of emerging artificial intelligence applications. Among the technologies considered promising for its realization, memristors have attracted considerable attention owing to their potential for enabling energy-efficient, brain-inspired information processing. A memristor is a two-terminal nonlinear electronic device that exhibits resistive switching and retains a memory of its previous states \cite{Mla19}. It mimics a biological synapse, whose conductance, corresponding to the synaptic weight, evolves in response to the activity of the connected neurons.

We consider oxide-based memristors, which typically consist of a thin transition-metal oxide layer (e.g., titanium dioxide) sandwiched between two metal electrodes \cite{Mla19}. In addition to conduction-band electrons and valence-band holes, these devices contain positively charged oxygen vacancies. When an electric field is applied, the vacancies migrate, leading to the formation of conductive filaments in oxygen-deficient regions of the oxide layer.

Memristor dynamics can be described at different modeling levels, ranging from simplified compact models commonly used in circuit simulations to atomistic and kinetic Monte Carlo approaches, which are computationally too demanding for device- or circuit-level simulations. The drift-diffusion equations offer a favorable compromise between computational tractability and physical accuracy. 

While the existence of weak solutions to memristor drift--diffusion systems has been established in the literature under various assumptions \cite{AGL25,JJZ23,JuVe24}, their qualitative properties have remained largely unexplored. In this paper, we derive time-uniform bounds and investigate the long-time behavior of solutions, showing convergence towards both equilibrium and non-equilibrium steady states.

\subsection{Model equations}

The (scaled) drift-diffusion equations for the electron density $n(x,t)$, hole density $p(x,t)$, oxygen-vacancy density $D(x,t)$, and electric potential $V(x,t)$ read as
\begin{align}
  \pa_t n &= \diver(\na n-n\na V) - R(n,p), \label{1.n} \\
  \pa_t p &= \diver(\na p+p\na V) - R(n,p), \label{1.p} \\
  \pa_t D &= \diver(\na D+D\na V), \label{1.D} \\
  \lambda^2\Delta V &= n-p-D+A(x)\quad\mbox{in }\Omega,\ t>0,
  \label{1.V}
\end{align}
where $\Omega\subset\R^d$ with $d=2$ is a bounded domain with Lipschitz boundary, $\lambda>0$ is the scaled Debye length, and $A(x)$ is the given density of immobile acceptor atoms. The equations are supplemented with the initial and mixed Dirichlet--Neumann boundary conditions
\begin{align}
  n(0) = n_I, \quad p(0) = p_I, \quad D(0) = D_I
  &\quad\mbox{in }\Omega, \label{1.ic} \\
  n=\bar{n}, \quad p=\bar{p}, \quad V=\bar{V} 
  &\quad\mbox{on }\Gamma_{\mathrm{D}}, \label{1.Dir} \\
  \na n\cdot\nu = \na p\cdot\nu = \na V\cdot\nu = 0
  &\quad\mbox{on }\Gamma_{\mathrm{N}}, \label{1.Neu} \\
  (\na D+D\na V)\cdot\nu = 0 &\quad\mbox{on }\pa\Omega. \label{1.bcD}
\end{align}
Here, $\pa\Omega=\Gamma_D\cup\Gamma_N$, $\Gamma_D$ models the Ohmic contacts and $\Gamma_N$ the insulating boundary parts. As the vacancies cannot leave the device, the no-flux boundary condition \eqref{1.bcD} implies the conservation of the total vacancy mass. The source term describes the generation and recombination of electron--hole pairs:
\begin{align*}
  R(n,p) = R_0(n,p)(np-1), 
\end{align*}
where $R_0$ is a positive, locally Lipschitz continuous function. If $R_0(n,p)=1/(1+\tau_p n + \tau_n p)$ for some capture-time constants $\tau_n$, $\tau_p>0$ and carrier concentrations $n$, $p$, the function $R$ is called the Shockley--Read--Hall term. 

\subsection{State of the art}

Although drift-diffusion equations for electrons and holes have been studied extensively in the literature \cite{Jue09,MRS90}, the mathematical analysis of multi-species models remains comparatively limited. For instance, global existence results in at most two space dimensions were established in \cite{BFPR14,GlHu05}. The dimensional restriction stems from the elliptic regularity $V\in W^{1,r}(\Omega)$ with $r>2$ \cite{Gro89}, which permits control of the drift term only in two dimensions. For the three-dimensional situation, the authors of \cite{ChLu95} assumed $W^{1,\infty}(\Omega)$ regularity for the Poisson equation, an assumption that can only be satisfied when the Dirichlet and Neumann portions of the boundary do not intersect. This result was subsequently improved in \cite{JJZ23}, where the weaker assumption $V\in W^{1,r}(\Omega)$ for some $r>2$ was sufficient to establish the existence of low-regularity solutions $n,p,D\in W^{1,1}(\Omega)$. The regularity of the solutions can be further improved by adopting high-density or Fermi--Dirac statistics \cite{AGL25,HJP25,JuVe24} in place of Maxwell--Boltzmann statistics, which lead to equations \eqref{1.n}--\eqref{1.D}.

The convergence of solutions to the two-species drift-diffusion system towards thermal equilibrium at an exponential rate was first investigated, to our knowledge, in \cite{GaGr89}. In the simpler setting of no-flux boundary conditions, the authors of \cite{FaIt95} established uniform-in-time estimates and constructed a compact global attractor. However, the relationship between this attractor and the stationary solutions remained unclear. Convergence to the thermal equilibrium, including a differential equation for the electrons on the trap level, was investigated in \cite{FeKn21}. More recently, \cite{KaSu22} established exponential convergence towards a time-periodic solution, thereby providing a long-time stability result for non-equilibrium steady states. For the two-species Nernst--Planck--Navier--Stokes system, one-dimensional steady states with nonzero currents are globally asymptotically stable if the currents are sufficiently small \cite{CIL22}. For related results for two-species Nernst--Planck--Stokes systems, see \cite{Lee25}. 

The analysis of the three-species system \eqref{1.n}--\eqref{1.V} is considerably more challenging than that of the two-species model for two main reasons. First, in the two-species setting, a monotonicity property can be exploited to reduce the quadratic drift term to a linear one. This property, however, does not extend to the three-species system. Second, the mismatch between the boundary conditions gives rise to boundary integrals, such as $\int_{\Gamma_D}D\na V\cdot\nu d\sigma$, that cannot be adequately controlled. In \cite{JJZ23}, these difficulties were overcome by establishing local $W^{1,q}(\omega)$ estimates for some $q>1$, where $\omega\Subset\Omega$. However, this level of regularity is insufficient to establish further qualitative properties of the solutions. This limitation can be overcome by restricting the analysis to two space dimensions. Indeed, in this setting, the regularity can be improved by exploiting the elliptic regularity $V\in W^{1,r}(\Omega)$ for some $r>2$, together with a logarithmic-type Gagliardo--Nirenberg inequality (see \eqref{1.logGN} below).

%%%%%%%%%%%%%%%%%

\subsection{Key ideas}\label{sec.ideas}

In this paper, we extend the results of \cite{JJZ23} by establishing time-uniform upper and lower bounds on the densities, as well as exponential convergence to both thermal equilibrium and non-equilibrium steady states. Our analysis relies on three new ideas.

First, the starting point for the a priori estimates is the energy dissipation inequality 
\begin{align}\label{1.dHdt}
  \frac{dH}{dt} &+ \int_\Omega\big(n|\na(\log n-V)|^2
  + p|\na(\log p+V)|^2 + D|\na(\log D+V)|^2\big)dx \\
  &\le C\Lambda(t)\bigg(1 + \int_\Omega(n+p)dx\bigg), \nonumber 
\end{align}
where $\Lambda(t)$, defined in \eqref{1.Lambda} below, depends on the time-dependent boundary data, $C>0$ is here and in the following a generic constant with values changing from line to line, and the (Helmholtz) free energy $H$ is the sum of the electron, hole, and vacancy entropies and the electric energy:
\begin{align}\label{1.energy}
  H(t) &= \int_\Omega\bigg(n\log\frac{n}{\bar{n}}-n+\bar{n}
  + p\log\frac{p}{\bar{p}}-p+\bar{p}
  + D\log\frac{D}{\bar{D}} - D + \bar{D}\bigg)dx \\
  &\phantom{xx}+ \frac{\lambda^2}{2}\int_\Omega
  |\na(V-\bar{V})|^2dx,
  \quad\mbox{where }\bar{D}:=\exp(-\bar{V}). \nonumber 
\end{align}
In the absence of the generation--recombination term, the total mass of electrons and holes is controlled, and hence the right-hand side of \eqref{1.dHdt} remains bounded. This is no longer the case in our setting. Indeed, estimating the integral $\int_\Omega(n+p)dx\le H+C_0$, for some constant $C_0>0$, would yield a bound on $H$ that depends on the length of the time interval $(0,T)$. To obtain a time-uniform estimate, our idea is to assume that $\int_0^\infty\Lambda(s)ds$ is finite. It follows from $\int_\Omega(n+p)dx\le H+C_0$ that $(d/dt)(H+C_0)\le C\Lambda(t)(H+C_0)$ and, by Gronwall's lemma,
\begin{align*}
  H(t)+C_0\le (H(0)+C_0)\exp\bigg(C\int_0^t\Lambda(s)ds\bigg)
  \le C\quad\mbox{for }t>0.
\end{align*}
This yields a priori estimates for $n,p,D$ in the Orlicz space $L\log L(\Omega)$, which are uniform in time if $\int_0^\infty\Lambda(t)dt<\infty$. From this point on, we can proceed as in \cite{JJZ23}. Indeed, we establish $L^2(\Omega)$ bounds for the densities using a logarithmic-type Gagliardo--Nirenberg inequality and then bootstrap these estimates to obtain $L^\infty(\Omega)$ bounds. The key Gagliardo--Nirenberg inequality is specific to two space dimensions, which explains our restriction $d=2$. 

We note that \cite[Theorem 1.2]{JJZ23} also claims time-uniform bounds. However, these bounds depend on the $L\log L(\Omega)$ norms of $n,p,D$, which are not, in general, uniformly bounded in time (because of the presence of the generation--recombination term).

Second, the lower bound for the vacancy density cannot be easily proved by a Stampacchia truncation argument as for $n,p$. Indeed, using $(D-m)_-=\min\{0,D-m\}$ for some $0<m\le\inf_\Omega D_I$ as a test function in the weak formulation of \eqref{1.D} yields, after an integration by parts, the boundary integral $\int_{\Gamma_D}\na V\cdot\nu(D-m)_-^2 d\sigma$, which cannot be controlled, since $D$ does not satisfy Dirichlet conditions on $\Gamma_D$. We overcome this issue by using, as in \cite[Lemma 3.6]{GaGr89}, the function $w:=-(\log(D+\delta)+m)_-\ge 0$ and the test function $-w^{q-1}/(D+\delta)$ with $\delta>0$, $q\ge 2$, $\exp(-m)\le\inf_\Omega D_I$, applying an Alikakos-type iteration argument, and passing to the limit $\delta\to 0$ to arrive at
\begin{align*}
  \sup_\Omega w(t)\le C_1:=C\max\Big\{1,\sup_{0<s<T} 
  \|w(s)\|_{L^1(\Omega)}\Big\} \quad\mbox{for }0<t<T.
\end{align*}
We infer that $\inf_\Omega D(t)\ge \exp(-C_1-m)>0$ if the $L^1(\Omega)$ norm of $w$ is bounded. Our new idea is to establish such a bound by estimating the Burg entropy $\int_\Omega(-\log D(t))dx$. 

Third, the long-time asymptotic behavior towards steady states is new. Even convergence to thermal equilibrium is nontrivial, since the time-dependent boundary conditions are generally not in thermal equilibrium and therefore require special treatment. An even more novel aspect is the analysis of the convergence towards {\em non-equilibrium} steady states. The idea is to exploit the assumption of small bias, measured by some parameter $\eta>0$. Indeed, the time derivative of the relative energy $H^*$, which measures the ``distance'' between $(n,p,D,V)$ and the steady state $(n^*,p^*,D^*,V^*)$, equals the negative dissipation $-\mathcal{D}^*$ (similarly as in \eqref{1.dHdt}) plus some cross terms, which can be bounded, leading to
\begin{align*}
  \frac{dH^*}{dt} + \mathcal{D^*} \le C\eta(\mathcal{D^*}+H^*).
\end{align*}
Together with the inequality $H^*\le C\mathcal{D}^*$, this yields $dH^*/dt + \kappa H^*\le 0$ with $\kappa=C^{-1}(1-C\eta)-C\eta$, which is positive for sufficiently small $\eta$, and hence exponential decay for $H^*(t)$. Since the constants depend on
the bounds of the solution, which are not known a priori (we do not have $\int_0^\infty\Lambda(t)dt<\infty$, since $\Lambda$ is constant for time-independent boundary data), we close the argument by a bootstrap argument. More precisely, $H^*$ controls the $L\log L$ norms (Lemma \ref{lem.HL}), these control the lower and upper bounds uniformly in time (Proposition \ref{prop.unifbd}), and therefore $H^*$ stays below the level $H^*(0) + 1$ for all times, providing the other bounds and closing the argument. 

\subsection{Main results}

We establish the following results. First, we prove the existence and uniqueness of a global weak solution to \eqref{1.n}--\eqref{1.bcD}. Compared with \cite{JJZ23}, our analysis accommodates generation--recombination terms and time-dependent boundary data, while extending the two-species framework of \cite{KaSu22} to a three-species system. Second, we establish uniform-in-time bounds for the densities and the electric energy. Third, we prove exponential convergence of the solution in two different regimes: (i) towards thermal equilibrium when the time-dependent boundary data converge to their equilibrium values as $t\to\infty$, and (ii) towards a non-equilibrium steady state with nonzero stationary fluxes when the boundary data are time-independent and the applied bias is sufficiently small.

Set $\Omega_T:=\Omega\times(0,T)$ for $T>0$ and $H_D^1(\Omega) := \{u\in H^1(\Omega):u=0\mbox{ on }\Gamma_D\}$. We define $z_+=\max\{0,z\}$, $z_-=\min\{0,z\}$, so that $z = z_+ + z_-$ for $z\in\R$. We impose the following assumptions:

\begin{itemize}
\item[(A1)] Domain: $\Omega\subset\R^d$ with $d=2$ is a bounded Lipschitz domain with $\pa\Omega=\Gamma_D\cup\Gamma_N$, where  $\Gamma_D\cap\Gamma_N = \emptyset$, $\Gamma_D$ has positive measure, and $\Gamma_N$ is relatively open in $\pa\Omega$. The boundary $\pa\Omega$ satisfies the regularity of \cite{Gro89}.
\item[(A2)] Data: $T>0$, $\lambda>0$, $A\in L^\infty(\Omega)$.
\item[(A3)] Generation--recombination term: $R(n,p)=R_0(n,p)(np-1)$, where $R_0\in C^0([0,\infty)^2)$ is positive, bounded, locally Lipschitz continuous, and satisfies $\inf_{[0,M]^2}R_0>0$ for every $M>0$. 
\item[(A4)] Initial data: $n_I,p_I,D_I\in L^\infty(\Omega)$ such that $n_I,p_I,D_I\ge m_I>0$ in $\Omega$. 
\item[(A5)] Boundary data: $\bar{n}=\exp(\bar{\mu}_n+\bar{V})$, $\bar{p}=\exp(\bar{\mu}_p-\bar{V})$, where $\Delta\bar{V}=0$ in $\Omega$, $\na\bar{V}\cdot\nu=0$ on $\Gamma_N$ weakly, and $|\bar{V}|\le\bar{M}$ in $\Omega\times(0,\infty)$ and, for $r>2$ given by \cite{Gro89},
\begin{align*}
  \sup_{t>0}\big(\|\bar{V}(t)\|_{W^{1,r}(\Omega)}
  + \|\pa_t\bar{V}(t)\|_{W^{1,r}(\Omega)}\big)
  + \|\bar{\mu}_n\|_{W^{1,\infty}(\Omega\times(0,\infty))}
  + \|\bar{\mu}_p\|_{W^{1,\infty}(\Omega\times(0,\infty))} < \infty.
\end{align*}
\end{itemize}

The existence of $r>2$ under Assumption (A1) is guaranteed by \cite{Gro89}: There exists $r>2$ such that the solution to $\Delta u=f$ in $\Omega\subset\R^2$, $u=\bar{V}$ on $\Gamma_D$, $\na u\cdot\nu = 0$ on $\Gamma_N$ satisfies $\|u\|_{W^{1,r}(\Omega)} \le C + C\|f\|_{L^{2r/(2+r)}(\Omega)}$. This result holds in two space dimensions only. We discuss extensions to three dimensions in Remark \ref{rem.3d}. Assumption (A3) includes the Shockley--Read--Hall term $R_0(n,p)=1/(1+\tau_pn+\tau_np)$. The positive lower bound for $R_0$ is not needed to derive finite-time bounds for the densities; it is only needed to prove time-uniform bounds; see Proposition \ref{prop.unifbd}. In Assumption (A5), the boundary data are parametrized by the harmonic extension $\bar{V}$ of the boundary potential and the quasi-Fermi potentials $\bar{\mu}_n=\log\bar{n}-\bar{V}$, $\bar{\mu}_p=\log\bar{p}+\bar{V}$. The data are in thermal equilibrium if $\bar{\mu}_n$, $\bar{\mu}_p$ are constant and $\bar{\mu}_n + \bar{\mu}_p = 0$. An applied bias corresponds to different values of the quasi-Fermi potentials on different contacts. The deviation of the boundary data is measured by 
\begin{align}\label{1.Lambda}
  \Lambda(t) &= \|(\na\bar{\mu}_n,\na\bar{\mu}_p)\|_{L^\infty(\Omega)}
  + \|(\pa_t\bar{\mu}_n,\pa_t\bar{\mu}_p)\|_{L^\infty(\Omega)}
  + \|\pa_t\bar{V}\|_{L^\infty(\Omega)} 
  + \|B(\bar{n}\bar{p})\|_{L^1(\Omega)},
\end{align}
where $B(b)=(b-1)(\log b)1_{\{b>1\}} + |\log b|1_{\{b\le 1\}}$ for $b>0$. Note that $\Lambda\equiv 0$ if and only if the boundary data are time-independent and in thermal equilibrium, i.e.\ $\na\bar{\mu}_n = \na\bar{\mu}_p = 0$ and $\bar{n}\bar{p}=1$.

\begin{theorem}[Global existence]\label{thm.ex}
Let Assumptions (A1)--(A5) hold. Then there exists a unique weak solution $(n,p,D,V)$ to \eqref{1.n}--\eqref{1.bcD} such that $n,p,D\ge 0$ in $\Omega_T$ and 
\begin{equation}\label{1.regul}
\begin{aligned}
  & n-\bar{n},p-\bar{p}\in L^2(0,T;H_D^1(\Omega)), \quad
  D\in L^2(0,T;H^1(\Omega)), \\
  & \pa_t n,\pa_t p\in L^2(0,T;H_D^1(\Omega)'), \quad
  \pa_t D\in L^2(0,T;H^1(\Omega)'), \\
  & V-\bar{V}\in L^\infty(0,T;H_D^1(\Omega)).
\end{aligned}
\end{equation}
Moreover, there exist constants $m_T$, $M_T>0$ such that $m_T\le n,p,D\le M_T$ in $\Omega_T$.  
\end{theorem}

The notion of weak solution is as usual. The solution is global in the sense that it exists on $(0,T)$ for arbitrary $T>0$. By uniqueness, the solutions on $(0,T)$ for different $T$ coincide on the common interval. As a consequence of the regularity stated in the theorem, $n,p,D\in C^0([0,T];L^2(\Omega))$, i.e., the initial data are satisfied in the sense of $L^2(\Omega)$. The theorem is proved by applying the Leray--Schauder fixed-point theorem to an approximate problem in which the nonlinearities are truncated. Compared to the proof of \cite{JJZ23}, we need to take care of the time-dependent boundary data and the generation--recombination term, which yields the free energy inequality \eqref{1.dHdt} without the dissipation integral. To obtain $L^2(\Omega)$ estimates, which are the starting point for the $L^\infty(\Omega)$ bounds, we need the restriction to two space dimensions. 

The upper bounds for the densities are derived by a Stampacchia truncation and Moser iteration argument. To obtain the lower bounds for $n$ and $p$, we only need a truncation argument and employ the Poisson equation in the drift term (no regularity for $\na V$ beyond $L^2(\Omega)$ is needed). As explained in Section \ref{sec.ideas}, this strategy cannot be applied to the vacancy density, and we need an Alikakos-type iteration together with an estimate for the Burg entropy and the regularity $\na V$ in $L^r(\Omega)$ with $r>2$. 

\begin{theorem}[Time-uniform bounds]\label{thm.bd}
Let Assumptions (A1)--(A5) hold and assume that $\Lambda\in L^1(0,\infty)$. Let $(n,p,D,V)$ be the weak solution to \eqref{1.n}--\eqref{1.bcD}. Then there exist constants $m,M>0$ such that
\begin{align*}
  0<m\le n,p,D\le M\quad\mbox{in }\Omega\times(0,\infty), \quad
  \sup_{t>0}\|V(t)\|_{W^{1,r}(\Omega)}<\infty.
\end{align*}
\end{theorem}

The assumption on $\Lambda$ is satisfied, for instance, if $\Lambda(t)\le C(1+t)^{-\alpha}$ for some $\alpha>1$. The lower and upper bounds established in Theorem \ref{thm.ex} depend on the length of the time interval. However, if $\Lambda\in L^1(0,\infty)$, the energy inequality \eqref{1.dHdt} yields time-uniform bounds for the densities. These bounds first imply uniform estimates in $L\log L(\Omega)$, which can then be successively improved to $L^2(\Omega)$ and, ultimately, $L^\infty(\Omega)$ bounds.

Finally, we prove two results on the long-time behavior of the solutions. First, we show the convergence to the thermal equilibrium state. A thermal equilibrium is a solution to \eqref{1.n}--\eqref{1.bcD} such that all fluxes and the generation--recombination term vanish. Assuming positive densities, this implies that $\na(\log n-V) = \na(\log p+V) = \na(\log D+V) = 0$ and $np = 1$. Since $\Omega$ is connected, there exist positive constants $c_n,c_p,c_D$ such that $\log n-V=\log c_n$, $\log p+V=\log c_p$, $\log D+V=\log c_D$. We set
\begin{align}\label{1.eq}
  n_{\rm eq} = c_n e^{V_{\rm eq}}, \quad 
  p_{\rm eq} = c_p e^{-V_{\rm eq}}, \quad 
  D_{\rm eq} = c_D e^{-V_{\rm eq}},
\end{align}
where $V_{\rm eq}\in H^1(\Omega)$ is the solution to
\begin{align}\label{1.Veq}
  & \lambda^2\Delta V_{\rm eq} = c_n e^{V_{\rm eq}}
  - c_p e^{-V_{\rm eq}} - c_D e^{-V_{\rm eq}} + A(x)
  \quad\mbox{in }\Omega, \\
  & V_{\rm eq} = \bar{V}_{\rm eq}\quad\mbox{on }\Gamma_D, \quad
  \na V_{\rm eq}\cdot\nu = 0 \quad\mbox{on }\Gamma_N. \label{1.Veqbc}
\end{align}
The constants $c_n,c_p$ (subject to $c_nc_p = n_{\rm eq}p_{\rm eq}=1$) and the boundary data $\bar{V}_{\rm eq}$ are given. In semiconductor theory, we have $c_n=c_p=1$ in the scaled equations, and $\bar{V}_{\rm eq}$ equals the built-in potential, which is determined from the doping profile. The constant $c_D$ is determined from
\begin{align}\label{1.cD}
  c_D = c_D[V_{\rm eq}] = \|D_I\|_{L^1(\Omega)}
  \bigg(\int_\Omega e^{-V_{\rm eq}}dx\bigg)^{-1},
\end{align}
which makes \eqref{1.Veq}--\eqref{1.Veqbc} an integro-differential problem. Therefore, the solution to \eqref{1.Veq}--\eqref{1.Veqbc} is nonstandard. We prove the existence of a weak solution $V_{\rm eq}$ in Proposition \ref{prop.Veq}. 

\begin{theorem}[Exponential convergence to thermal equilibrium]
\label{thm.exp1}
Let Assumptions (A1)--(A5) hold and assume that 
\begin{align}\label{1.A6} 
  \Lambda(t) + \|\pa_t\bar{V}(t)\|_{H^1(\Omega)}\le Ce^{-\kappa_0 t},
  \quad t>0,
\end{align}
for some $C>0$ and $\kappa_0>0$. Let $(n,p,D,V)$ be the weak solution to \eqref{1.n}--\eqref{1.bcD} and $(n_{\rm eq},p_{\rm eq},D_{\rm eq},V_{\rm eq})$ be the weak solution to \eqref{1.eq}--\eqref{1.cD} with $\log c_n:=\lim_{t\to\infty}\bar\mu_n(t)$, $c_p:=1/c_n$, and $\bar{V}_{\rm eq}:=\lim_{t\to\infty}\bar{V}(t)$ (these limits exist by \eqref{1.A6}). Then there exist constants $C>0$ and $0<\kappa_1<\kappa_0$ such that for all $t>0$,
\begin{align*}
  \|n(t)-n_{\rm eq}\|_{L^2(\Omega)}
  + \|p(t)-p_{\rm eq}\|_{L^2(\Omega)}
  + \|D(t)-D_{\rm eq}\|_{L^2(\Omega)} 
  + \|V(t)-V_{\rm eq}\|_{H^1(\Omega)}\le Ce^{-\kappa_1 t}.
\end{align*}
\end{theorem}

Condition \eqref{1.A6} is crucial. If this expression is converging algebraically, then we obtain only algebraic convergence towards thermal equilibrium. Theorem \ref{thm.exp1} is proved by the relative entropy method. The relative energy is defined by
\begin{align*}
 H^* = \int_\Omega\big(h(n|n^*) + h(p|p^*) + h(D|D^*)\big)dx
  + \frac{\lambda^2}{2}\int_\Omega|\na(V-V^*)|^2 dx,
\end{align*}
where $h(u|u^*) = u(\log(u/u^*)-1)+u^*$ and $(n^*,p^*,D^*,V^*)$ is the steady state (here the thermal equilibrium). Differentiating $H^*$ with respect to time, we find that $dH^*/dt + \mathcal{D}^*\le C\mathcal{S}(t)$, where $\mathcal{D}^*$ is the dissipation integral and $\mathcal{S}(t)$ contains the boundary data and decays like $\exp(-\kappa_0 t)$ by assumption. Showing that $H^*\le C\mathcal{D}^*$ and applying Gronwall's lemma concludes the proof. The technical difficulty is that the thermal equilibrium does not satisfy the Dirichlet data, which requires the introduction of boundary-adapted equilibrium densities using a harmonic extension operator.

The final result concerns the long-time asymptotics towards a non-equilibrium steady state. We say that $(n^*,p^*,D^*,V^*)\in H^1(\Omega)^4$ is a steady state if the stationary versions of \eqref{1.n}--\eqref{1.V}, the boundary conditions \eqref{1.Dir}--\eqref{1.bcD}, and $m^*\le n^*,p^*,D^*\le M^*$ hold in $\Omega$ for some $m^*,M^*>0$. We show in Section \ref{sec.exp2} that the steady state is a solution to the stationary bipolar drift--diffusion equations with an additional nonlocal term in the Poisson equation, similarly as \eqref{1.cD}. 

\begin{theorem}[Exponential convergence to non-equilibrium steady state]
\label{thm.exp2}
Let Assumptions (A1)--(A5) hold with time-independent boundary data. Let $r>2$ as in Assumption (A5). Let $(n,p,D,V)$ be the weak solution to \eqref{1.n}--\eqref{1.bcD} and $(n^*,p^*,D^*,V^*)$ be a steady state with the same boundary data. Then there exist ${\eta_0}\in(0,1]$, $\kappa_1>0$, $C>0$, such that if
\begin{align*}
  \|\na(\log n^*-V^*)\|_{L^r(\Omega)}
  + \|\na(\log p^*+V^*)\|_{L^r(\Omega)}
  + \|n^*p^*-1\|_{L^\infty(\Omega)} \le {\eta_0},
\end{align*}
then
\begin{align*}
  \|n(t)-n^*\|_{L^2(\Omega)} + \|p(t)-p^*\|_{L^2(\Omega)}
  + \|D(t)-D^*\|_{L^2(\Omega)} + \|V(t)-V^*\|_{H^1(\Omega)}
  \le Ce^{-\kappa_1 t}.
\end{align*}
The value ${\eta_0}$ depends on $H^*(0)$, $m^*$, $M^*$, and the data.
\end{theorem}

The parameter $\eta$ measures the distance of the steady state from thermal equilibrium. In thermal equilibrium, we have $\eta=0$. For a small applied bias, we expect $\eta$ to be small, but we do not prove an estimate for $\eta$ in terms of the boundary data here. 

The theorem is proved by a variant of the relative entropy method. The relative energy inequality yields a right-hand side involving $\na\mu_n^*$, $\na\mu_p^*$ and a remainder generation--recombination term. These expressions can be bounded by $C\eta(H^*+\mathcal{D}^*)$, leading to $dH^*/dt + c H^*\le C\eta H^*$, which yields exponential decay for $H^*$ if $\eta>0$ is chosen sufficiently small. However, the bounds are not known a priori, since $\Lambda\not\in L^1(0,\infty)$ for non-equilibrium data, so Theorem \ref{thm.bd} does not apply. As long as $H^*(t)\le C(H^*(0)+1)$, our bounds hold with fixed constants. Consequently, taking into account the relative energy inequality, $H^*$ is decreasing. A continuity argument shows that this situation persists for all times, and we conclude the argument. 

The paper is organized as follows. Section \ref{sec.aux} collects some auxiliary results. The existence and uniqueness result, Theorem \ref{thm.ex}, is proved in Sections \ref{sec.ex} and \ref{sec.bd}. Section \ref{sec.unifbd} is devoted to the proof of the time-uniform bounds stated in Theorem \ref{thm.bd}. Finally, the asymptotic stability results, Theorems \ref{thm.exp1} and \ref{thm.exp2}, are established in Sections \ref{sec.exp1} and \ref{sec.exp2}, respectively. 

%%%%%%%%%%%%%%%%%%%%%%%%%%%%%%%%%%%%%%%%%%%%%%%%%%%%%%%%%%%%%%

\section{Auxiliary results}\label{sec.aux}

We recall some auxiliary lemmas used in this paper.

\begin{lemma}[Gagliardo--Nirenberg inequality \cite{Zei90}]
Let $\Omega\subset\R^2$ be a bounded domain with Lipschitz boundary. Then there exists $C>0$ such that for all $u\in H^1(\Omega)$,
\begin{align}\label{1.GN}
  \|u\|_{L^q(\Omega)} &\le C\|u\|_{H^1(\Omega)}^{1-2/q}
  \|u\|_{L^2(\Omega)}^{2/q} \quad\mbox{for } 2\le q<\infty, \\
  \|u\|_{L^q(\Omega)} &\le C\|u\|_{H^1(\Omega)}^{1-1/q}
  \|u\|_{L^1(\Omega)}^{1/q} \quad\mbox{for } 1\le q<\infty.
  \label{1.GNa}
\end{align}
\end{lemma}

\begin{lemma}[Logarithmic-type Gagliardo--Nirenberg inequality \cite{BHN94,GlHu97}]
Let $\Omega\subset\R^2$ be a bounded domain with Lipschitz boundary, $q>1$, and $\delta>0$. Then there exists $C(\delta)>0$ such that for all $u\in H^1(\Omega)$,
\begin{align}\label{1.logGN}
  \|u\|_{L^q(\Omega)}^q \le \delta\|u\|_{H^1(\Omega)}^{q-1}
  \|u\log|u|\|_{L^1(\Omega)} + C(\delta)\|u\|_{L^1(\Omega)}. 
\end{align}
\end{lemma}

\begin{lemma}[Uniform $L^\infty(\Omega)$ bounds \cite{HoJu20}]
\label{lem.ali}
Let $w$ be a nonnegative function such that $w^{q/2}\in L^2(0,T;H^1(\Omega))\cap C^0([0,T];L^2(\Omega))$ for all $q\ge 2$. Assume that there exist $K\ge 1$ and $\gamma\ge 0$ such that, for all $q\ge 2$ and $0<t<T$,
\begin{align*}
  e^t\int_\Omega w(t)^q dx 
  + \int_0^t e^s\int_\Omega |\na w^{q/2}|^2 dxds
  \le Kq^\gamma\bigg(e^t + \int_0^t e^s\int_\Omega w^q dxds\bigg).
\end{align*}
Then there exists a constant $C(K,\gamma)>0$ not depending on $T$ such that
\begin{align*}
  \sup_{0<t<T}\|w(t)\|_{L^\infty(\Omega)}
  \le C(K,\gamma)\max\Big\{1,\sup_{0<t<T}\|w(t)\|_{L^1(\Omega)}\Big\}.
\end{align*}
In particular, if $w\in L^\infty(0,\infty;L^1(\Omega))$ and the assumption holds with $T=\infty$, then also $w\in L^\infty(\Omega\times(0,\infty))$. 
\end{lemma}

\begin{lemma}\label{lem.aux}
Let $a,b>0$. Then
\begin{align*}
  (1-a)(\log a-\log b) \le B(b) := (b-1)\log b1_{\{b>1\}}
  + |\log b|1_{\{b\le 1\}}.
\end{align*}
\end{lemma}

\begin{proof}
If $a\le b\le 1$ or $1\le b\le a$ or $a\le 1\le b$ or $b\le 1\le a$, we have $(1-a)(\log a-\log b)\le 0$. Let $b\le a\le 1$. Then $(1-a)(\log a-\log b)\le -(1-a)\log b\le -\log b$. Finally, if $1\le a\le b$, it follows that $(1-a)(\log a-\log b)=(a-1)(\log b-\log a)\le (a-1)\log b\le (b-1)\log b$.
\end{proof}

%%%%%%%%%%%%%%%%%%%%%%%%%%%%%%%%%%%%%%%%%%%%%%%%%%%%%%%%%%%%%%

\section{Existence of a global bounded weak solution}\label{sec.ex}

We prove the existence of a weak solution on finite time intervals. For this, we solve first an approximate problem with truncation parameter $k\in\N$, derive uniform estimates from the free energy \eqref{1.energy}, and pass to the limit $k\to\infty$. 

\subsection{Solution of an approximate problem}

We truncate the drift and generation--recombination terms by introducing the cutoff functions $T_k(u)=\max\{0,\min\{k,u\}\}$ for $k>0$ and $u\in\R$. Set $R_k(n,p)=R(T_k(n),T_k(p))$. The approximate problem reads as
\begin{align}
  \pa_t n &= \diver(\na n-T_k(n)\na V) - R_k(n,p), \label{2.n} \\
  \pa_t p &= \diver(\na p+T_k(p)\na V) - R_k(n,p), \label{2.p} \\
  \pa_t D &= \diver(\na D+T_k(D)\na V), \label{2.D} \\
  \lambda^2\Delta V &= n-p-D+A(x)\quad\mbox{in }\Omega,\ t>0,
  \label{2.V}
\end{align}
together with the initial and boundary conditions \eqref{1.ic}--\eqref{1.bcD}. We assume throughout this section that Assumptions (A1)--(A5) hold. 

\begin{lemma}[Global existence for the approximate problem]
%\label{lem.approx}
There exists a weak solution $(n_k,p_k,D_k)$ to \eqref{1.ic}--\eqref{1.bcD}, \eqref{2.n}--\eqref{2.V} such that $n_k,p_k,D_k\ge 0$ in $\Omega_T$ and 
the regularity \eqref{1.regul} holds.
\end{lemma}

\begin{proof}
The result is proved similarly as in \cite[Sec.\ 2.2]{JJZ23}, but we need to take care of the time-dependent boundary conditions and the source terms. We only sketch the proof, highlighting the treatment of the source terms. Let $(\widetilde n,\widetilde p,\widetilde D)\in L^2(\Omega_T;\R^3)$ and $\sigma\in[0,1]$. The linear system
\begin{align}
  \lambda^2\Delta V &= \widetilde n - \widetilde p - \widetilde D 
  + A(x), \nonumber \\ %\label{2.VV} \\
  \pa_t n &= \diver(\na n-\sigma T_k(\widetilde n)\na V) 
  - \sigma R_k(\widetilde n,\widetilde p), 
  \label{2.nn} \\
  \pa_t p &= \diver(\na p+\sigma T_k(\widetilde p)\na V) 
  - \sigma R_k(\widetilde n,\widetilde p), 
  \label{2.pp} \\
  \pa_t D &= \diver(\na D+\sigma T_k(\widetilde D)\na V)
  \quad\mbox{in }\Omega,\ t>0, \label{2.D2} 
\end{align}
with the initial and boundary conditions 
\begin{align*}
  n(0) = \sigma n_I, \quad p(0) = \sigma p_I, \quad D(0) = \sigma D_I
  &\quad\mbox{in }\Omega, \\
  n=\sigma\bar{n}, \quad p=\sigma\bar{p}, \quad V=\bar{V} 
  &\quad\mbox{on }\Gamma_{\mathrm{D}}, \\
  \na n\cdot\nu = \na p\cdot\nu = \na V\cdot\nu = 0
  &\quad\mbox{on }\Gamma_{\mathrm{N}}, \\
  (\na D+\sigma T_k(\widetilde D)\na V)\cdot\nu = 0 
  &\quad\mbox{on }\pa\Omega
\end{align*}
admits a unique weak solution $(n,p,D,V)\in L^2(0,T;H^1(\Omega))^4$, since the drift term satisfies $|T_k(\widetilde n)\na V|\in L^2(\Omega_T)$ and $R_k$ is bounded, so standard linear parabolic theory can be applied \cite{Zei90}. This defines the fixed-point operator $F:L^2(\Omega_T)^3\times[0,1] \to L^2(\Omega_T)^3$, $(\widetilde n,\widetilde p,\widetilde D;\sigma)\mapsto(n,p,D)$. The map is continuous, satisfies $F(\widetilde n,\widetilde p,\widetilde D,0)=0$, and is, because of standard energy-type estimates and the Aubin--Lions lemma, compact. It remains to derive a uniform bound for all fixed points of $F(\cdot;\sigma)$. We use the test functions $n-\sigma\bar{n}$ in the weak formulation of \eqref{2.nn} with $(\widetilde n,\widetilde p)=(n,p)$. Since the drift term can be estimated as $|\int_\Omega T_k(n)\na V\cdot\na ndx|\le k\|\na V\|_{L^2(\Omega)}\|\na n\|_{L^2(\Omega)}$, we obtain, using Young's inequality to absorb $\na n$ into the left-hand side and using the standard elliptic estimate for $\na V$ (see \cite[Sec.\ 2.2]{JJZ23}), 
\begin{align*}
  \frac12&\int_\Omega (n-\sigma\bar{n})(t)^2 dx
  + \int_0^t\int_\Omega|\na n|^2 dxds
  \le C(k) + C(k)\int_0^t\int_\Omega(n^2+p^2+D^2)dxds \\
  &+ \frac{\sigma}{2}\int_\Omega(n_I-\bar{n}(0))^2 dx
  - \sigma\int_0^t\int_\Omega\pa_t\bar{n}(n-\sigma\bar{n})dxds
  - \sigma\int_0^t\int_\Omega R_k(n,p)(n-\sigma\bar{n})dxds.
\end{align*}
The last terms on the right-hand side were not present in \cite[Sec.\ 2.2]{JJZ23}; they can be estimated in a standard way, using Young's inequality and the boundedness of $R_k$. Estimating in a similar way for the equations for $p$ and $D$, we end up with
\begin{align*}
  \int_\Omega(n^2 &+ p^2 + D^2)(t)dx
  + \int_0^t\int_\Omega\big(|\na n|^2 + |\na p|^2 + |\na D|^2\big)dxds \\
  &\le C(k) + C(k)\int_0^t\int_\Omega(n^2 + p^2 + D^2)dxds.
\end{align*}
Gronwall's inequality yields a bound for $n,p,D$ in $L^\infty(0,T;L^2(\Omega))$ uniform in $\sigma$. This provides the uniform bound required for the Leray--Schauder fixed-point theorem. This theorem yields the existence of a fixed point of $F(\cdot,1)$, which is a solution to \eqref{2.n}--\eqref{2.V}.

It remains to prove the nonnegativity of the densities. We take the test function $n_- = \min\{0,n\}$ in the weak formulation of \eqref{2.n} and use $n(0)_-=0$, $\bar{n}_-=0$, as well as $T_k(n)=0$ and $R_k(n,p)=-R_0(0,T_k(p))\le 0$ on $\{n<0\}$:
\begin{align*}
  \frac12\int_\Omega n_-(t)^2 dx + \int_0^t\int_\Omega|\na n_-|^2 dxds
  = \int_0^t\int_\Omega R_0(0,T_k(p))n_- dxds \le 0.
\end{align*}
This shows that $n\ge 0$ in $\Omega_T$. 
\end{proof}

%%%%%%%%%%%%%%%%%%%%

\subsection{Uniform estimates}

Next, we derive some estimates uniform in $k$. Let $0\le\delta<1$. As in \cite[Sec.\ 2.3]{JJZ23}, we introduce the functions
\begin{align*}
  g_{k,\delta}(u) = \int_1^u\int_1^y\frac{dzdy}{T_k(z)+\delta}, 
  \quad 
  G_{k,\delta}(u|\bar{u}) = g_{k,\delta}(u) - g_{k,\delta}(\bar{u})
  - g'_{k,\delta}(\bar{u})(u-\bar{u}) 
\end{align*}
for $u\in\R$, $\bar{u}>0$ and the regularized free energy
\begin{align*}
  H_{k,\delta} = \int_\Omega\bigg(
  G_{k,\delta}(n|\bar{n}) + G_{k,\delta}(p|\bar{p}) 
  + G_{k,\delta}(D|\bar{D}) + \frac{\lambda^2}{2}
  |\na(V-\bar{V})|^2\bigg)dx,
\end{align*} 
where $\bar{D}:=\bar{D}_{k,\delta}$ is uniquely defined by $g'_{k,\delta}(\bar{D})=-\bar{V}$. The function $g'_{k,\delta}(u)$ approximates $\log u$, and the number $\bar{D}$ is uniformly bounded with respect to $k$ and $\delta$. The function $g'_{k,\delta}$ is constructed such that $(T_k(z)+\delta)g_{k,\delta}''(z) = 1$ for $z\in\R$. 
%We set $H_k:=H_{k,0}$ and $h_k:=h_{k,0}$.

\begin{lemma}[Regularized free energy inequality]\label{lem.H}
Let $H_k:=H_{k,0}$ and $k\ge k_0 = \max\{1,$ $\|\bar n\|_{L^\infty(\Omega)}, \|\bar p\|_{L^\infty(\Omega)}, \|n_I\|_{L^\infty(\Omega)}, \|p_I\|_{L^\infty(\Omega)}, \|D_I\|_{L^\infty(\Omega)}\}$. Let $(n_k,p_k,D_k,V_k)$ be a weak solution to \eqref{2.n}--\eqref{2.V}, \eqref{1.ic}--\eqref{1.bcD}. Then there exists a constant $C>0$ independent of $k$ such that
\begin{align}\label{2.Hk}
  H_{k}(t) \le H(0) + C\int_0^t\bigg(\Lambda(s)+\Lambda^2(s)+\frac{1}{k}\bigg)
  \bigg(1 + \int_\Omega(n_k+p_k)dx\bigg)ds,
\end{align}
where $\Lambda(t)$ is defined in \eqref{1.Lambda}. Consequently, $H_k(t)\le C(T)$ for $t\in(0,T)$, where $C(T)$ is independent of $k$ but depends on $T$. 
\end{lemma}

\begin{proof}
Set $\phi_n=g'_{k,\delta}(n_k)-V_k$, $\phi_p=g'_{k,\delta}(p_k)+V_k$, $\phi_D=g'_{k,\delta}(D_k)+V_k$ and $\bar{\phi}_n = g'_{k,\delta}(\bar{n})-\bar{V}$, $\bar{\phi}_p = g'_{k,\delta}(\bar{p})+\bar{V}$. A computation shows that
\begin{align}\label{2.aux3}
  \frac{dH_{k,\delta}}{dt} &= \langle\pa_t n_k,\phi_n-\bar{\phi}_n\rangle
  + \langle\pa_t p_k,\phi_p-\bar{\phi}_p\rangle
  + \langle\pa_t D_k,\phi_D\rangle \\
  &\phantom{xx}- \int_\Omega\big(
  (n_k-\bar{n})\pa_t g_{k,\delta}'(\bar{n})
  + (p_k-\bar{p})\pa_t g_{k,\delta}'(\bar{p}) 
  - (D_k-\bar{D})\pa_t\bar{V}\big)dx. \nonumber 
\end{align}
We insert equation \eqref{2.n} for $n_k$ and use Young's inequality:
\begin{align*}
  \langle\pa_t n_k,\phi_n-\bar{\phi}_n\rangle
  &= -\int_\Omega\big((T_k(n_k)+\delta)\na\phi_n+\delta\na V_k\big)
  \cdot\na(\phi_n-\bar{\phi}_n)dx \\ 
  &\phantom{xx}- \int_\Omega R_k(n_k,p_k)(\phi_n-\bar{\phi}_n)dx \\
  &\le -\frac12\int_\Omega(T_k(n_k)+\delta)
  |\na(\phi_n-\bar{\phi}_n)|^2 dx
  + \int_\Omega(T_k(n_k)+\delta)|\na\bar{\phi}_n|^2 dx \\
  &\phantom{xx}+ \delta\int_\Omega|\na V_k|^2 dx 
  - \int_\Omega R_k(n_k,p_k)(\phi_n-\bar{\phi}_n)dx.
\end{align*}
The same estimate holds for $p_k$, whereas $\langle\pa_t D_k,\phi_D\rangle \le -\frac12\int_\Omega(T_k(D_k)+\delta)|\na\phi_D|^2 dx + \frac{\delta}{2}\int_\Omega|\na V_k|^2 dx$. The last integral in \eqref{2.aux3} is bounded by $C\Lambda(1+\int_\Omega(n_k+p_k)dx)$, since $\bar{n}$ and $\bar p$ are strictly positive. 

Next, we consider the source term. Combining the terms in the equations for $n_k$ and $p_k$, we need to estimate 
\begin{align}
  X &:= (1-T_k(n_k)T_k(p_k))
  \big(g'_{k,\delta}(n_k)+g'_{k,\delta}(p_k)
  -g'_{k,\delta}(\bar{n})-g'_{k,\delta}(\bar{p})\big) \nonumber\\
  &= (1-a)(\log a_\delta - \log\bar{a}_\delta)
  + (1-a)\frac{(n_k-k)_+ + (p_k-k)_+}{k+\delta}, \label{X}
\end{align}
where $a:=T_k(n_k)T_k(p_k)$,  $a_\delta := (T_k(n_k)+\delta)(T_k(p_k)+\delta)$, and $\bar{a}_\delta := (\bar{n}+\delta)(\bar{p}+\delta)$. Let $k\ge\max\{1, \|\bar n\|_{L^\infty(\Omega)}, \|\bar p\|_{L^\infty(\Omega)}\}$, so that $T_k(\bar n)=\bar n$ and $T_k(\bar p)=\bar p$. We claim that  
\begin{align}\label{2.X}
  X \le B(\bar{a}_\delta) 
  + \bigg(\delta+\frac{1}{k}\bigg)(n_k+p_k) + \delta^2,
\end{align}
where $B$ is defined in Lemma \ref{lem.aux}. Indeed, if $a>1$, we have $\log a_\delta>\log a$ and the first and second term in \eqref{X} are each nonpositive. Thus, $X\le 0$. If $a\le 1$, we use 
\begin{align*}
  (1-a)\log a_\delta \le (a_\delta-1)_+ 
  \le \delta(n_k+p_k)+\delta^2, \quad
  -(1-a)\log\bar{a}_\delta
  \le -(\log\bar{a}_\delta)_-\le B(\bar{a}_\delta)
\end{align*}
to infer directly \eqref{2.X}, proving the claim. This shows in particular that
\begin{align*}
  -&\int_\Omega R_k(n_k,p_k)\big(g'_{k,\delta}(n_k)+g'_{k,\delta}(p_k)
  -g'_{k,\delta}(\bar{n})-g'_{k,\delta}(\bar{p})\big)dx \\
  &= \int_\Omega R_0(T_k(n),T_k(p))X dx
  \le C\bigg(\frac{1}{k} + \delta\bigg)\bigg(1+\int_\Omega(n_k+p_k)dx\bigg)
  +  C\int_\Omega B(\bar{a}_\delta)dx.
\end{align*}

Summarizing, we have proved that 
\begin{align}\label{2.aux4}
  H_{k,\delta}&(t) - H_{k,\delta}(0) \le \int_0^t\int_\Omega
  \big((T_k(n_k)+\delta)|\na\bar{\phi}_n|^2
  + (T_k(p_k)+\delta)|\na\bar{\phi}_p|^2\big)dxds \\
  &\phantom{xx}+ \delta\int_0^t\int_\Omega\big(3|\na V_k|^2
  + |\na\bar{\phi}_n|^2 + |\na\bar{\phi}_p|^2\big)dxds \nonumber \\
  &\phantom{xx}+ C\int_0^t\bigg(\Lambda + \frac{1}{k} + \delta\bigg)
  \bigg(1 + \int_\Omega(n_k+p_k)dx\bigg)ds
  + C\int_0^t\int_\Omega|B(\bar{a}_\delta)|dxds. \nonumber 
\end{align}
We pass to the limit $\delta\to 0$ in this inequality, keeping $k$ fixed. Using $|G_{k,\delta}(u|\bar{u})|\le C(k)(1+u^2)$ for $0<c\le\bar{u}\le C$, dominated convergence implies that $H_{k,\delta}(t)\to H_{k}(t)$ for $t\ge 0$. Since 
\begin{align*}
  \na\bar{\phi}_n = \frac{\bar{n}}{\bar{n}+\delta}\na
  (\log\bar{n}-\bar{V}) - \frac{\delta}{\bar{n}+\delta}\na\bar{V}
  \to \na(\log\bar{n}-\bar{V})\quad\mbox{a.e. as }\delta\to 0
\end{align*}
and $T_k(n)\le k$, $T_k(\bar{n})=\bar{n}$, we can again apply dominated convergence to find that 
\begin{align*}
  \int_0^t\int_\Omega(T_k(n_k)+\delta)|\na\bar{\phi}_n|^2 dxds
  \to \int_0^t\int_\Omega T_k(n_k)|\na(\log\bar{n}-\bar{V})|^2 dxds
  \le C\int_0^t\int_\Omega\Lambda^2 n_kdxds.
\end{align*}
The limit $\delta\to 0$ in \eqref{2.aux4} then proves \eqref{2.Hk}. 
\end{proof}

The free energy inequality \eqref{2.Hk} and the inequality $G_{k,0}(u|\bar{u})\ge \frac12 u\log u-C$ yield the following estimates uniform in $k$.

\begin{lemma}[Uniform $L\log L$ bounds]\label{lem.LlogL}
Let $(n_k,p_k,D_k,V_k)$ be a weak solution to \eqref{2.n}--\eqref{2.V}, \eqref{1.ic}--\eqref{1.bcD}. Then there exists $C(T)>0$ independent of $k$ such that
\begin{align*}
  \sup_{0<t<T}\int_\Omega\big(
  |n_k\log n_k| + |p_k\log p_k| + |D_k\log D_k|
  + |\na V_k|^2\big)(t)dx \le C(T).
\end{align*}
In particular, $n_k$ and $p_k$ are uniformly bounded in $L^\infty(0,T;L^1(\Omega))$. 
\end{lemma}

Up to now, the results hold for any space dimension. To obtain uniform $L^2(\Omega)$ bounds, we need the restriction $d=2$, which comes from the logarithmic-type Gagliardo--Nirenberg inequality \eqref{1.logGN}. 

\begin{lemma}[Uniform $L^2(\Omega)$ bounds]\label{lem.L2H1}
Let $(n_k,p_k,D_k,V_k)$ be a weak solution to \eqref{2.n}--\eqref{2.V}, \eqref{1.ic}--\eqref{1.bcD}. Then there exists $C(T)>0$ independent of $k$ such that
\begin{align*}
  \|(n_k,p_k,D_k)\|_{L^\infty(0,T;L^2(\Omega))}
  + \|(n_k,p_k,D_k)\|_{L^2(0,T;H^1(\Omega))}
  + \|V_k\|_{L^\infty(0,T;W^{1,r}(\Omega))} \le C(T),
  %\|(\pa_tn_k,\pa_tp_k)\|_{L^2(0,T;H^1_D(\Omega)')}
  %+ \|\pa_tD_k\|_{L^2(0,T;H^1(\Omega)')}
  %&\le C(T),
\end{align*}
recalling that $r>2$ by Assumption (A5).
\end{lemma}

\begin{proof}
We use the test function $n_k-\bar{n}$ in the weak formulation of \eqref{2.n}:
\begin{align*}
  \frac12&\frac{d}{dt}\int_\Omega(n_k-\bar{n})^2 dx
  + \int_\Omega|\na(n_k-\bar{n})|^2 dx
  = -\int_\Omega\pa_t\bar{n}(n_k-\bar{n})dx \\
  &- \int_\Omega\na\bar{n}\cdot\na(n_k-\bar{n})dx
  + \int_\Omega T_k(n_k)\na V_k\cdot\na(n_k-\bar{n})dx
  - \int_\Omega R_k(n_k,p_k)(n_k-\bar{n})dx.
\end{align*}
By the bounds from Lemma \ref{lem.LlogL}, the first term on the right-hand side is bounded by $C+C\|n_k-\bar{n}\|_{L^2(\Omega)}^2$, using Young's inequality. For the second term, again by Young's inequality, we absorb the $|\na(n_k-\bar{n})|^2$ term into the left-hand side. For the reaction term, we estimate $-R_k(n_k,p_k)(n_k-\bar{n})\le Cn_k$ if $n_k>\bar{n}$, while for $n_k\le\bar{n}\le k$, we have
\begin{align*}
  -R_k(n_k,p_k)(n_k-\bar{n}) = R_0(T_k(n_k),T_k(p_k))(n_kT_k(p_k)-1)(\bar{n}-n_k)
  \le C(n_k+p_k).
\end{align*}
Thus, in both cases, $-\int_\Omega R_k(n_k,p_k)(n_k-\bar n)dx\le C(T)$, by Lemma \ref{lem.LlogL}. It remains to estimate the drift term. We use H\"older's inequality with $1/q+1/2+1/r=1$ (i.e.\ $q=2r/(r-2)$):
\begin{align}\label{2.drift}
  \bigg|\int_\Omega T_k(n_k)\na V_k\cdot\na(n_k-\bar{n})dx\bigg|
  &\le C\|n_k-\bar{n}\|_{L^q(\Omega)}\|\na V_k\|_{L^r(\Omega)}
  \|\na(n_k-\bar{n})\|_{L^2(\Omega)} \\
  &\phantom{xx}+ C\|\bar{n}\|_{L^q(\Omega)}\|\na V_k\|_{L^r(\Omega)}
  \|\na(n_k-\bar{n})\|_{L^2(\Omega)}. \nonumber 
\end{align}
We apply the logarithmic-type Gagliardo--Nirenberg inequality \eqref{1.logGN} and then the Poincar\'e inequality: 
\begin{align*}
  \|n_k-\bar{n}\|_{L^q(\Omega)}
  &\le \|n_k\|_{L^q(\Omega)} + \|\bar{n}\|_{L^q(\Omega)} \\
  &\le \delta C\|n_k\|_{H^1(\Omega)}^{1-1/q}
  \|n_k\log n_k\|_{L^1(\Omega)}^{1/q} 
  + C(\delta)\|n_k\|_{L^1(\Omega)}^{1/q} + \|\bar{n}\|_{L^q(\Omega)} \\
  &\le C\delta \|\na(n_k-\bar{n})\|_{L^2(\Omega)}^{(r+2)/(2r)} 
  + C(\delta).
\end{align*}
By the elliptic regularity $V\in W^{1,r}(\Omega)$ and the Gagliardo--Nirenberg inequality \eqref{1.GNa},
\begin{align}\label{2.Vaux}
  \|\na V_k\|_{L^r(\Omega)} 
  &\le C + C\|n_k-p_k-D_k\|_{L^{2r/(2+r)}(\Omega)} \\
  &\le C + C\|n_k-p_k-D_k\|_{H^1(\Omega)}^\theta
  \|n_k-p_k-D_k\|_{L^1(\Omega)}^{1-\theta} \nonumber \\
  &\le C + C\|\na(n_k-p_k-D_k)\|_{L^2(\Omega)}^{(r-2)/(2r)}.
  \nonumber 
\end{align}
These estimations allow us to treat the first term on the right-hand side of \eqref{2.drift}:
\begin{align*}
  \|n_k&-\bar{n}\|_{L^q(\Omega)}\|\na V_k\|_{L^r(\Omega)}
  \|\na(n_k-\bar{n})\|_{L^2(\Omega)} \\
  &\le C\delta\big(\|\na(n_k-\bar{n})\|_{L^2(\Omega)}^{2}
  + \|\na(p_k-\bar{p})\|_{L^2(\Omega)}^{2} 
  + \|\na D_k\|_{L^2(\Omega)}^{2}\big) + C(\delta).
\end{align*}
In three space dimensions, the exponent is larger than two, which is the reason why we need to suppose that $d=2$. The last term in \eqref{2.drift} is estimated similarly:
\begin{align*}
  \|\bar{n}\|_{L^q(\Omega)}\|\na V_k\|_{L^r(\Omega)}
  \|\na(n_k-\bar{n})\|_{L^2(\Omega)}
  &\le \delta\big(\|\na(n_k-\bar{n})\|_{L^2(\Omega)}^2
  + \|\na(p_k-\bar{p})\|_{L^2(\Omega)}^2 \\
  &\phantom{xx}+ \|\na D_k\|_{L^2(\Omega)}^2\big) + C(\delta).
\end{align*}
We conclude from \eqref{2.drift} that
\begin{align*}
  \bigg|\int_\Omega & T_k(n_k)\na V_k\cdot\na(n_k-\bar{n})dx\bigg| \\
  &\le C\delta\big(\|\na(n_k-\bar{n})\|_{L^2(\Omega)}^{2}
  + \|\na(p_k-\bar{p})\|_{L^2(\Omega)}^{2} 
  + \|\na D_k\|_{L^2(\Omega)}^{2}\big) + C(\delta).
\end{align*}
Similar estimates are derived for $p_k$ and $D_k$. Thus, for sufficiently small $\delta>0$, 
\begin{align*}
  \frac{d}{dt}\int_\Omega&
  \big((n_k-\bar n)^2+(p_k-\bar p)^2+D_k^2\big)dx
  + \int_\Omega(|\na n_k|^2 + |\na p_k|^2 + |\na D_k|^2)dx \\
  &\le C(\delta) + C\int_\Omega(n_k+p_k)dx \le C,
\end{align*}
which finishes the proof.
\end{proof}

\begin{lemma}[Further uniform bounds]\label{lem.further}
Let $(n_k,p_k,D_k,V_k)$ be a weak solution to \eqref{2.n}--\eqref{2.V}, \eqref{1.ic}--\eqref{1.bcD}. Then there exists $C(T)>0$ independent of $k$ such that
\begin{align*}
  %\|\na V_k\|_{L^\infty(0,T;L^r(\Omega))}
  %+ \|n_kp_k\|_{L^2(\Omega_T)} &\le C(T), \\
  \|(\pa_t n_k,\pa_t p_k)\|_{L^2(0,T;H_D^1(\Omega)')}
  + \|\pa_t D_k\|_{L^2(0,T;H^1(\Omega)')} &\le C(T).
\end{align*}
\end{lemma}

\begin{proof}
The Gagliardo--Nirenberg inequality implies that $\|n_k\|_{L^4(\Omega)}\le C\|n_k\|_{H^1(\Omega)}^{1/2} \|n_k\|_{L^2(\Omega)}^{1/2}$, so $n_k$ and $p_k$ are uniformly bounded in $L^4(\Omega_T)$, hence $n_kp_k$ and also $R_k(n_k,p_k)$ are bounded in $L^2(\Omega_T)$ uniformly in $k$. Furthermore, with $q=2r/(r-2)$,
\begin{align*}
  \|T_k(n_k)\na V_k\|_{L^2(\Omega_T)}
  &\le \|n_k\|_{L^2(0,T;L^{q}(\Omega))}
  \|\na V_k\|_{L^\infty(0,T;L^r(\Omega))} \\
  &\le C\|n_k\|_{L^2(0,T;H^1(\Omega))}
  \|\na V_k\|_{L^\infty(0,T;L^r(\Omega))} \le C(T),
\end{align*}
since $H^1(\Omega)\hookrightarrow L^q(\Omega)$ in two space dimensions. The bounds for the time derivative follow from \eqref{2.n}--\eqref{2.D}.
\end{proof}

%%%%%%%%%%%%%%%%%

\subsection{Limit $k\to\infty$}

Lemmas \ref{lem.L2H1} and \ref{lem.further} allow us to apply the Aubin--Lions lemma to conclude the existence of a subsequence (not relabeled) such that 
\begin{align*}
  (n_k,p_k,D_k)\to (n,p,D) &\quad\mbox{strongly in }L^2(\Omega_T)
  \mbox{ and a.e.\ in }\Omega_T, \\
  (n_k,p_k,D_k)\rightharpoonup (n,p,D) &\quad\mbox{weakly in }
  L^2(0,T;H^1(\Omega)), \\
  (\pa_tn_k,\pa_tp_k,\pa_tD_k)\rightharpoonup
  (\pa_t n,\pa_t p,\pa_t D) &\quad\mbox{weakly in }
  L^2(0,T;H_D^1(\Omega)')^2\times L^2(0,T;H^1(\Omega)').
\end{align*}
Taking into account the Poisson equation, this implies that $V_k\to V$ strongly in $L^2(0,T;$ $H^1(\Omega))$. Since $|T_k(n_k)-n|\le|T_k(n_k)-T_k(n)| + |T_k(n)-n|\le |n_k-n|+(n-k)_+$, we have $T_k(n_k)\to n$ strongly in $L^2(\Omega_T)$. Together with the strong convergence of $\na V_k$, this gives $T_k(n_k)\na V_k\to n\na V$ strongly in $L^1(\Omega_T)$, and since the sequence is bounded in $L^2(\Omega_T)$, also weakly in $L^2(\Omega_T)$. Similar results hold for $p_k$ and $D_k$. Furthermore, $R_k(n_k,p_k)\to R(n,p)$ a.e.\ in $\Omega_T$ (by the continuity of $R_0$), and $(R_k(n_k,p_k))$ is bounded in $L^2(\Omega_T)$, so $R_k(n_k,p_k)\rightharpoonup R(n,p)$ weakly in $L^2(\Omega_T)$. Passing to the limit in \eqref{2.n}--\eqref{2.V}, we conclude that $(n,p,D,V)$ is a weak solution to \eqref{1.n}--\eqref{1.bcD} on $(0,T)$.

%%%%%%%%%%%%%%%%%%%%%%%%%%%%%%%%%%%%%%%%%%%%%%%%%%%%%%%%%%%%

\section{Boundedness and positivity}\label{sec.bd}

We show that every weak solution is bounded and strictly positive on finite time intervals. 

\subsection{Upper and lower bounds}

We first derive upper bounds for the densities. They are needed to show the lower bounds. 

\begin{lemma}\label{lem.up}
Let $(n,p,D,V)$ be a weak solution to \eqref{1.n}--\eqref{1.bcD} on $(0,T)$. Then there exists $M>0$ such that
\begin{align*}
  \sup_{0<t<T}\big(\|n(t)\|_{L^\infty(\Omega)}
  + \|p(t)\|_{L^\infty(\Omega)} + \|D(t)\|_{L^\infty(\Omega)}\big)
  \le M,
\end{align*}
where $M$ depends on $T$ and the $L^\infty(0,T;L^2(\Omega))$ norm of $(n,p,D)$. 
\end{lemma}

\begin{proof}
We use an Alikakos-type iteration together with a truncation argument. Let $M$ be larger than the $L^\infty(\Omega)$ norms of $n_I,p_I,D_I, \bar{n},\bar{p}$, and let $K>M$, $q\ge 2$. We set 
\begin{align*}
  w = (\min\{n,K\}-M)_+, \quad
  \Phi(n) = \int_0^n(\min\{z,K\}-M)_+^{q-1}dz,
\end{align*}
recalling that $z_+=\max\{0,z\}$, and use $w^{q-1}=\Phi'(n)$ as a test function in the weak formulation of \eqref{1.n} (note that $\Phi(n(0))=0$, since $n_I\le M$):
\begin{align}\label{2.aux5}
  \frac{d}{dt}\int_\Omega\Phi(n)dx
  + \frac{4(q-1)}{q^2}\int_\Omega|\na w^{q/2}|^2 dx
  = \int_\Omega n\na V\cdot\na w^{q-1}dx - \int_\Omega R(n,p)w^{q-1}dx.
\end{align}
The last term is estimated according to $-R(n,p)w^{q-1}\le Cw^{q-1}\le C(1+w^q)$. For the drift term, we compute
\begin{align*}
  \bigg|\int_\Omega n\na V\cdot\na w^{q-1}dx\bigg|
  &= \bigg|\int_{\{M<n<K\}}(w+M)\na V\cdot\na w^{q-1}dx\bigg| \\
  &= \frac{2(q-1)}{q}\bigg|\int_\Omega(w^{q/2}+Mw^{q/2-1})
  \na V\cdot\na w^{q/2}dx\bigg| \nonumber \\
  &\le 2\int_\Omega\big(w^{q/2} + M(1+w^{q/2})\big)
  |\na V||\na w^{q/2}|dx. \nonumber 
\end{align*}
By H\"older's inequality with $1/\sigma+1/r+1/2=1$ (i.e., $\sigma=2r/(r-2)$), the Gagliardo--Nirenberg inequality \eqref{1.GN} with $\theta = 1-2/\sigma=2/r$, and Young's inequality with exponents $2/(1+\theta)$, $2/(1-\theta)$ and $\delta>0$,
\begin{align*}
  \int_\Omega &w^{q/2}|\na V||\na w^{q/2}|dx
  \le \|\na V\|_{L^r(\Omega)}\|w^{q/2}\|_{L^{\sigma}(\Omega)}
  \|\na w^{q/2}\|_{L^2(\Omega)} \\
  &\le C(V)\|w^{q/2}\|_{H^1(\Omega)}^{1+\theta}
  \|w^{q/2}\|_{L^2(\Omega)}^{1-\theta}
  \le \delta\|\na w^{q/2}\|_{L^2(\Omega)}^2 
  + C(\delta,V)\|w^{q/2}\|_{L^2(\Omega)}^2, \\
  \int_\Omega &|\na V||\na w^{q/2}|dx
  \le \|\na V\|_{L^2(\Omega)}\|\na w^{q/2}\|_{L^2(\Omega)}
  \le \delta\|\na w^{q/2}\|_{L^2(\Omega)}^2 + C(V)\delta^{-\gamma}
\end{align*}
for some $\gamma>0$, where $C(V)$ depends on the $L^\infty(0,T;W^{1,r}(\Omega))$ norm of $V$ and hence on the $L^\infty(0,T;L^2(\Omega))$ norm of $(n,p,D)$. The $\delta$-terms can be absorbed into the left-hand side of \eqref{2.aux5} if we choose $\delta= c/q$ for some sufficiently small $c>0$. Integration over $(0,t)$ and using the inequality $w^q\le q\Phi(n)$, we deduce from \eqref{2.aux5} that
\begin{align*}
  \int_\Omega w(t)^q dx + \int_0^t\int_\Omega|\na w^{q/2}|^2 dxds
  \le Cq^{\gamma+1}\bigg(\int_0^t\int_\Omega w^q dxds + t\bigg).
\end{align*}
We multiply this inequality by $e^t$:
\begin{align*}
  e^t\int_\Omega w(t)^q dx 
  + \int_0^t e^s \int_\Omega|\na w^{q/2}|^2 dxds
  \le Cq^{\gamma+1}
  \bigg(e^T\int_0^t e^s\int_\Omega w^q dx + Te^t\bigg).
\end{align*}
We deduce from Lemma \ref{lem.ali} the bound $\|w(t)\|_{L^\infty(\Omega)}\le C_T$ for $0<t<T$, where $C_T>0$ depends on the $L^1(\Omega)$ norm of $n$ (and hence on $T$), but not on $K$. The limit $K\to\infty$ implies that $(n-M)_+\le C_T$ a.e., giving the upper bound. The proof for $p$ and $D$ is similar.
\end{proof}

We prove next that the electron and hole densities are strictly positive.

\begin{lemma}[Lower bounds for $n$ and $p$]\label{lem.lownp}
Let $(n,p,D,V)$ be a weak solution to \eqref{1.n}--\eqref{1.bcD} on $(0,T)$. Then there exists $m_T>0$, depending on $T$ and the upper bound $M$ from Lemma \ref{lem.up}, such that $n\ge m_T$ and $p\ge m_T$ in $\Omega_T$. 
\end{lemma}

\begin{proof}
Let $m_0=\min\{\inf_\Omega n_I,\inf_{\Omega\times(0,T)}\bar n\}>0$, $\mu>0$ and set $m(t)=m_0\exp(-\mu t)$. We use the test function $(n-m)_-\le 0$ in the weak formulation of \eqref{1.n}, taking into account that $\pa_t m=-\mu m$ and $(n-m)_-=0$ on $\Gamma_D$,
\begin{align*}
  \frac12&\frac{d}{dt}\int_\Omega(n-m)_-^2 dx
  + \int_\Omega|\na(n-m)_-|^2 dx 
  = \int_\Omega(n-m)\na V\cdot\na(n-m)_- dx \\
  &+ m\int_\Omega\na V\cdot\na(n-m)_- dx - \int_\Omega R(n,p)(n-m)_- dx
  + \mu m\int_\Omega(n-m)_- dx.
\end{align*}
Using the Poisson equation \eqref{1.V}, the first two terms on the right-hand side can be written as
\begin{align*}
  \int_\Omega(n&-m)\na V\cdot\na(n-m)_- dx
  + m\int_\Omega\na V\cdot\na(n-m)_- dx \\
  &= \frac12\int_\Omega\na V\cdot\na(n-m)_-^2 dx 
  + m\int_\Omega\na V\cdot\na(n-m)_- dx \\
  &= -\frac{1}{2\lambda^2}\int_\Omega(n-p-D+A)(n-m)_-(n+m) dx \\
  &\le -\frac{m}{\lambda^2}(m+\|A\|_{L^\infty(\Omega)})
  \int_\Omega(n-m)_- dx \le -Cm\int_\Omega(n-m)_- dx,
\end{align*}
where $C>0$ also depends on $\|A\|_{L^\infty(\Omega)}$. For the source term, we use the assumption that $R_0$ is bounded and $p\le M$ by Lemma \ref{lem.up}. Then
\begin{align*}
  -\int_\Omega R(n,p)(n-m)_-dx 
  = \int_{\{n<m\}}R_0(n,p)(np-1)|(n-m)_-|
  \le -CMm\int_\Omega(n-m)_- dx
\end{align*}
We conclude that
\begin{align*}
  \frac12\frac{d}{dt}\int_\Omega(n-m)_-^2 dx
  + \int_\Omega|\na(n-m)_-|^2 dx 
  \le \big(C(1+M)-\mu\big)m\int_\Omega|(n-m)_-|dx
  \le 0,
\end{align*}
choosing $\mu\ge C(1+M)$. This implies that $(n-m)_-=0$ and $n\ge m$ in $\Omega_T$.
\end{proof}

The vacancy density requires a different argument. Indeed, we cannot integrate by parts in the drift term using the Poisson equation, because the test function $D$ does not vanish on $\Gamma_D$ and the corresponding boundary integral is not controlled. Instead, we estimate the drift term via H\"older's inequality and not via the Poisson equation. First, we derive an inequality for the so-called Burg entropy $\int_\Omega(-\log D)dx$ \cite{Bur72}.

\begin{lemma}[Burg entropy bound for $D$]\label{lem.burg}
Let $(n,p,D,V)$ be a weak solution to \eqref{1.n}--\eqref{1.bcD}. Then 
\begin{align}\label{3.burg}
  \int_\Omega(-\log D(t))dx 
  + \frac12\int_0^t\int_\Omega|\na\log D|^2 dxds
  \le \int_\Omega(-\log D_I)dx 
  + \frac12\int_0^t\|\na V\|_{L^2(\Omega)}^2ds.
\end{align}
\end{lemma} 

\begin{proof}
Let $\delta>0$. We use the test function $-1/(D+\delta)$ in \eqref{1.D}:
\begin{align*}
  \frac{d}{dt}\int_\Omega(-\log(D+\delta))dx
  &= -\int_\Omega(\na D+D\na V)\cdot\frac{\na D}{(D+\delta)^2}dx \\
  &\le -\|\na\log(D+\delta)\|_{L^2(\Omega)}^2
  + \|\na V\|_{L^2(\Omega)}\|\na\log(D+\delta)\|_{L^2(\Omega)}. 
\end{align*}
Young's inequality and integration over $(0,t)$ show \eqref{3.burg} for $D+\delta$ instead of $D$. Since $-\log(D+\delta)$ is monotone in $\delta$, monotone convergence and Fatou's lemma allow us to pass to the limit $\delta\to 0$.
\end{proof}

The Poisson equation and the regularity $n,p,D\in L^\infty(0,T;L^2(\Omega))$ provide a bound for $\na V$ in $L^\infty(0,T;L^2(\Omega))$ so that $\int_\Omega(-\log D(t))dx$ is bounded on $(0,T)$, implying that $D(t)>0$ a.e.\ in $\Omega_T$. We prove that $D(t)$ is even bounded from below by a positive constant.

\begin{lemma}[Lower bound for $D$]\label{lem.lowD}
Let $(n,p,D,V)$ be a weak solution to \eqref{1.n}--\eqref{1.bcD}. Then $D\ge m_T$ in $\Omega_T$, where $m_T>0$ depends on $\sup_{0<t<T} \int_\Omega(-\log D(t))dx$ and hence on $T$.
\end{lemma}

\begin{proof}
Let $w:=-(\log(D+\delta)+m)_-\ge 0$. We use the test function $-w^{q-1}/(D+\delta)$ in the weak formulation of \eqref{1.D}, where $\delta>0$ and $m>0$ is such that $\exp(-m)\le \inf_\Omega D_I$. Let $q\ge 2$. Then
\begin{align*}
  \frac{1}{q}\frac{d}{dt}&\int_\Omega w^qdx 
  = -\int_\Omega (\na D+D\na V)\cdot\frac{\na D}{(D+\delta)^2}
  \big((q-1)w^{q-2}+w^{q-1}\big)dx \\
  &= -\int_\Omega\bigg(|\na w|^2 - \frac{D}{D+\delta}
  \na V\cdot\na w\bigg)\big((q-1)w^{q-2}+w^{q-1}\big)dx \\
  &= -\frac{4(q-1)}{q^2}\int_\Omega |\na w^{q/2}|^2 dx
  - \frac{4}{(q+1)^2}\int_\Omega |\na w^{(q+1)/2}|^2 dx \\
  &\phantom{xx}+ \frac{2}{q}\int_\Omega \frac{D}{D+\delta}
  \big((q-1)w^{q/2-1} + w^{q/2}\big)\na w^{q/2}\cdot\na V dx.
\end{align*}
Taking into account that $w^{q/2-1}\le 1+w^{q/2}$, we obtain
\begin{align}\label{2.wq}
  \frac{1}{q}&\frac{d}{dt}\int_\Omega w^qdx
  + \frac{4(q-1)}{q^2}\int_\Omega |\na w^{q/2}|^2 dx
  \le 2\int_\Omega (1+w^{q/2})|\na w^{q/2}||\na V|dx \\
  &\le 2\|\na w^{q/2}\|_{L^2(\Omega)}
  \|\na V\|_{L^2(\Omega)} + 2\int_\Omega w^{q/2}|\na w^{q/2}|
  |\na V| dx =: I_1+I_2. \nonumber 
\end{align}
Since $n,p,D$ are bounded in $L^2(\Omega)$ uniformly in $(0,T)$, also $\|\na V\|_{L^2(\Omega)}$ is bounded uniformly in $(0,T)$, so that, by Young's inequality, $I_1 \le \frac{1}{2q}\|\na w^{q/2}\|_{L^2(\Omega)}^2 + Cq$. For the second term, we use \eqref{2.Vaux}, giving $\|\na V\|_{L^r(\Omega)}\le C$. We use H\"older's inequality with $1/r+1/\sigma+1/2=1$ (thus $\sigma=2r/(r-2)$), Gagliardo--Nirenberg's inequality with $\theta=2/r<1$, and Young's inequality with exponents $2/(1+\theta)$ and $2/(1-\theta)$:
\begin{align*}
  I_2 &\le 2\|w^{q/2}\|_{L^{\sigma}(\Omega)}
  \|\na w^{q/2}\|_{L^2(\Omega)}\|\na V\|_{L^r(\Omega)}
  \le C\|w^{q/2}\|_{H^1(\Omega)}^\theta
  \|w^{q/2}\|_{L^2(\Omega)}^{1-\theta}\|\na w^{q/2}\|_{L^2(\Omega)} \\
  &\le \frac{1}{2q}\|\na w^{q/2}\|_{L^2(\Omega)}^2
  + Cq^{\gamma}\|w^{q/2}\|_{L^2(\Omega)}^2
\end{align*}
for some $\gamma>0$. We infer from \eqref{2.wq} after integration over $t$ that
\begin{align*}
  \int_\Omega w(t)^q dx 
  + \int_0^t \int_\Omega|\na w^{q/2}|^2 dxds
  \le Cq^2t + Cq^{1+\gamma}\int_0^t\int_\Omega w^q dxds.
\end{align*}
Thus, we can proceed as in the proof of Lemma \ref{lem.up} (using $w(0)=0$, since $\exp(-m)\le\inf_\Omega D_I$) to conclude that $\sup_\Omega w(t)\le C_1:=C\max\{1,\sup_{0<t<T} \|w(t)\|_{L^1(\Omega)}\}$. The limit $\delta\to 0$ leads to $-\inf_\Omega(\log D(t)+m)_-\le C_1$ and hence to $\inf_\Omega D(t)\ge \exp(-C_1-m)>0$ for $0<t<T$. The constant $C_1$ depends on $\|(\log D(t)+m)_-\|_{L^1(\Omega)}$. It follows from $-(\log y)_- = -\log y+(\log y)_+\le -\log y+y$ that
\begin{align*}
  \|(\log D(t)+m)_-\|_{L^1(\Omega)}
  \le |\Omega||m| + \int_\Omega(-\log D(t))dx + \|D_I\|_{L^1(\Omega)},
\end{align*}
concluding the proof.
\end{proof}

%%%%%%%%%%%%%%%%%%%%

\subsection{Uniqueness of a weak solution}

We prove that the bounded weak solution to \eqref{1.n}--\eqref{1.bcD} is unique. Let $(n_i,p_i,D_i,V_i)$ for $i=1,2$ be two weak solutions to \eqref{1.n}--\eqref{1.bcD} on $(0,T)$. We know that the densities are bounded by $M$ and that $\|\na V_i\|_{L^\infty(0,T;L^r(\Omega))}\le C$, where $r>2$. Let $n:=n_1-n_2$, $p:=p_1-p_2$, $D:=D_1-D_2$, and $V:=V_1-V_2$. We use $n,p,D$ as test functions in the weak formulation of \eqref{1.n}--\eqref{1.D}, respectively, and sum the equations:
\begin{align}\label{2.uniq}
  \frac12\frac{d}{dt}&\int_\Omega(n^2+p^2+D^2)dx
  + \int_\Omega(|\na n|^2 + |\na p|^2 + |\na D|^2)dx \\
  &= \int_\Omega(n\na V_1+n_2\na V)\cdot\na n dx
  - \int_\Omega(p\na V_1+p_2\na V)\cdot\na p dx \nonumber \\
  &\phantom{xx}- \int_\Omega(D\na V_1+D_2\na V)\cdot\na D dx
  - \int_\Omega(R(n_1,p_1)-R(n_2,p_2))(n+p)dx. \nonumber 
\end{align}
We estimate the first drift term by using H\"older's inequality with $1/q+1/r+1/2=1$ and then the Gagliardo--Nirenberg inequality \eqref{1.GN} with $\theta=2/r<1$:
\begin{align*}
  \bigg|\int_\Omega n\na V_1\cdot\na n dx\bigg|
  &\le \|n\|_{L^q(\Omega)}\|\na V_1\|_{L^r(\Omega)}
  \|\na n\|_{L^2(\Omega)} \\
  &\le C\|n\|_{H^1(\Omega)}^{1+\theta}\|n\|_{L^2(\Omega)}^{1-\theta}
  \le \frac14\|\na n\|_{L^2(\Omega)}^2 + C\|n\|_{L^2(\Omega)}^2, \\
  \bigg|\int_\Omega n_2\na V\cdot\na n dx\bigg|
  &\le \|n_2\|_{L^\infty(\Omega)}\|\na V\|_{L^2(\Omega)}
  \|\na n\|_{L^2(\Omega)} \\
  &\le \frac14\|\na n\|_{L^2(\Omega)}^2 
  + C(M)\big(\|n\|_{L^2(\Omega)}^2 + \|p\|_{L^2(\Omega)}^2
  + \|D\|_{L^2(\Omega)}^2\big).
\end{align*}
We deduce from \eqref{2.uniq} and the Lipschitz continuity of $R$ on $[0,M]^2$ that
\begin{align*}
  \frac12\frac{d}{dt}&\int_\Omega(n^2+p^2+D^2)dx
  + \frac12\int_\Omega(|\na n|^2 + |\na p|^2 + |\na D|^2)dx \\
  &\le C\big(\|n\|_{L^2(\Omega)}^2 + \|p\|_{L^2(\Omega)}^2
  + \|D\|_{L^2(\Omega)}^2\big).
\end{align*}
Since $n(0)=p(0)=D(0)=0$, Gronwall's lemma implies that $n=p=D=0$ and then $V=0$ in $\Omega_T$, showing that the weak solution is indeed unique.

%%%%%%%%%%%%%%%%%%%%%%%%%%%%%%%%%%%%%%%%%%%%%%%%%%%%%%%%%%%%

\section{Time-uniform bounds}\label{sec.unifbd}

In this section, we prove time-uniform upper and lower bounds for the densities under the assumption $\Lambda\in L^1(0,\infty)$. We introduce the quasi-Fermi potentials
\begin{align*}
  \mu_n = \log n-V, \quad \mu_p = \log p+V, \quad \mu_D = \log D+V,
\end{align*}
as well as $\bar{\mu}_n=\log\bar{n}-\bar{V}$, $\bar{\mu}_p = \log\bar{p}+\bar{V}$. Thanks to the lower bounds in Lemmas \ref{lem.lownp} and \ref{lem.lowD} and Assumption (A5), they are well defined. 

\begin{lemma}[Free energy inequality]\label{lem.fei}
Let $(n,p,D,V)$ be a weak solution to \eqref{1.n}--\eqref{1.bcD} on $(0,T)$. Then there exists a constant $C>0$, not depending on $T$, such that
\begin{align}\label{3.fei}
  \frac{dH}{dt} + \frac12\int_\Omega\big(n|\na\mu_n|^2 + p|\na\mu_p|^2
  + D|\na\mu_D|^2\big)dx
  \le C\Lambda(t)\bigg(1 + \int_\Omega(n+p)dx\bigg),
\end{align}
where $\Lambda(t)$ is defined in \eqref{1.Lambda}.
\end{lemma}

Compared to Lemma \ref{lem.H}, we provide dissipation bounds. This is possible, since we can use the lower bounds for the densities to define the quasi-Fermi potentials.

\begin{proof}
We use the test functions $\mu_n-\bar{\mu}_n$, $\mu_p-\bar{\mu}_p$, $\mu_D$ in the weak formulation of \eqref{1.n}--\eqref{1.D}, respectively. Then, similarly as in \eqref{2.aux3},
\begin{align}\label{3.dHdt} 
  \frac{dH}{dt} &= \langle\pa_t n,\mu_n-\bar{\mu}_n\rangle
  + \langle\pa_t p,\mu_p-\bar{\mu}_p\rangle 
  + \langle\pa_t D,\mu_D\rangle \\
  &\phantom{xx}- \int_\Omega\big((n-\bar{n})\pa_t\log\bar{n}
  + (p-\bar{p})\pa_t\log\bar{p} - (D-\bar{D})\pa_t\bar{V}\big)dx. \nonumber 
\end{align}
We infer from $\pa_t n=\diver(n\na\mu_n)-R(n,p)$ and Young's inequality that
\begin{align*}
  \langle\pa_t n,\mu_n-\bar{\mu}_n\rangle
  &= -\int_\Omega n|\na\mu_n|^2 dx
  + \int_\Omega n\na\mu_n\cdot\na\bar{\mu}_n dx
  - \int_\Omega R(n,p)(\mu_n-\bar{\mu}_n)dx \\
  &\le -\frac12\int_\Omega n|\na\mu_n|^2 dx 
  + C\Lambda\int_\Omega ndx
  - \int_\Omega R(n,p)(\mu_n-\bar{\mu}_n)dx,
\end{align*}
and similarly for $p$ (the computation for $D$ is even simpler). The source terms combine to
\begin{align*}
  -\int_\Omega &R(n,p)(\mu_n-\bar{\mu}_n + \mu_p-\bar{\mu}_p)dx \\
  &= \int_\Omega R_0(n,p)(1-np)\big(\log(np)-\log(\bar{n}\bar{p})\big)dx
  \le C\|B(\bar{n}\bar{p})\|_{L^1(\Omega)} \le C\Lambda,
\end{align*}
where $B$ is defined in Lemma \ref{lem.aux}. Finally, we have $|\pa_t\log\bar{n}|\le |\pa_t\bar{\mu}_n|+|\pa_t\bar{V}|\le\Lambda$, so the last integral in \eqref{3.dHdt} is bounded by $C\Lambda(1+\int_\Omega(n+p)dx)$. 
\end{proof}

The following result follows from definition \eqref{1.energy} of $H$.

\begin{lemma}[Consequences of the free energy bound]\label{lem.cons}
Let $(n,p,D,V)$ be a weak solution to \eqref{1.n}--\eqref{1.bcD} on $(0,T)$. Then there exists $C_*\ge 1$, only depending on the data, such that $H(t)+C_*\ge 1$, $\int_\Omega(n+p)(t)dx\le H(t)+C_*$, and
\begin{align}\label{5.L}
  L(t) := \int_\Omega\big(|n\log n| + |p\log p| + |D\log D|
  + |\na V|^2\big)(t)dx \le C_*(H(t)+C_*).
\end{align}
\end{lemma}

\begin{proof}
We infer from $\int_\Omega D\bar{V}dx\ge -C\|D_I\|_{L^1(\Omega)}$ that the entropy parts and $\|\na(V-\bar{V})\|_{L^2(\Omega)}^2$ are bounded by $H+C$. The inequality $u\le h(u|\bar{u}) + (e-1)\bar{u}$ gives the bound for $\int_\Omega(n+p)dx$. We estimate
\begin{align*}
  |n\log n| \le n\log n + 2/e 
  = h(n|\bar{n}) + n\log\bar{n}+n-\bar{n}+2/e
  \le h(n|\bar{n}) + n(|\log\bar{n}|+1)
\end{align*}
and similarly for $p$ and $D$. The result follows from $\|\na V\|_{L^2(\Omega)} \le |\na(V-\bar{V})\|_{L^2(\Omega)}+\|\na\bar{V}\|_{L^2(\Omega)}$.
\end{proof}

\begin{proposition}[Free energy bound]\label{prop.feb}
Let $(n,p,D,V)$ be a weak solution to \eqref{1.n}--\eqref{1.bcD} on $(0,T)$. Then, with the constants $C$ from Lemma \ref{lem.fei} and $C_*$ from Lemma \ref{lem.cons},
\begin{align*}
  H(t) + C_* \le (H(0)+C_*)\exp\bigg(2C\int_0^t\Lambda(s)ds\bigg)
  \quad\mbox{for }0<t<T.
\end{align*}
In particular, if $T=\infty$ and $\Lambda\in L^1(0,\infty)$, then $H(t)$ is uniformly bounded for $t\in(0,\infty)$. 
\end{proposition}

\begin{proof}
It follows from Lemma \ref{lem.cons} that $1+\int_\Omega(n+p)dx\le 2(H+C_*)$. Inserting this bound into \eqref{3.fei} shows that $\frac{d}{dt}(H+C_*)\le 2C\Lambda(t)(H+C_*)$, and an integration completes the proof.
\end{proof}

We prove now some bounds only depending on $L(t)$ from Lemma \ref{lem.cons}. Thus, they are uniform in time if $\sup_{t>0}L(t)<\infty$ or if $\sup_{t>0}H(t)<\infty$ or if $\Lambda\in L^1(0,\infty)$.

\begin{lemma}[Time-uniform $L^2(\Omega)$ bound]%\label{lem.unifL2}
Let $(n,p,D,V)$ be a weak solution to \eqref{1.n}--\eqref{1.bcD} on $(0,T)$ and assume that $L(t)\le L_*$ for $0<t<T$. Then
\begin{align*}
  \|n(t)\|_{L^2(\Omega)}^2 + \|p(t)\|_{L^2(\Omega)}^2
  + \|D(t)\|_{L^2(\Omega)}^2 \le C(L_*),
\end{align*}
where $C(L_*)>0$ depends only on $L_*$ and the data, but not on $T$.
\end{lemma}

\begin{proof}
Similarly as in the proof of Lemma \ref{lem.L2H1}, we use $n-\bar{n}$ as a test function in the weak formulation of \eqref{1.n}:
\begin{align*}
  \frac12\frac{d}{dt}\int_\Omega(n-\bar{n})^2dx
  &+ \int_\Omega|\na n|^2 dx
  = \int_\Omega\na n\cdot\na\bar{n}dx
  - \int_\Omega\pa_t\bar{n}(n-\bar{n})dx \\
  &+ \int_\Omega n\na V\cdot\na(n-\bar{n})dx
  - \int_\Omega R(n,p)(n-\bar{n})dx.
\end{align*}
The first two terms on the right-hand side can be bounded by $C(L_*)$ in a standard way, using the gradient term on the left-hand side. For the source term, if $n\ge\bar{n}$, we have $-R(n,p)(n-\bar{n})=R_0(n,p)(1-np)(n-\bar{n})\le C(n-\bar{n})$, while for $n\le\bar{n}$, it holds that $-R(n,p)(n-\bar{n})=R_0(n,p)(1-np)(n-\bar{n}) \le Cnp\bar{n}\le C\bar{n}^2 p$. We conclude that
\begin{align*}
  -\int_\Omega R(n,p)(n-\bar{n}) \le C\int_\Omega(n+p)dx\le C(L_*). 
\end{align*}
The drift term is estimated as in the proof of Lemma \ref{lem.L2H1} (using $d=2$), leading to
\begin{align*}
  \bigg|\int_\Omega n\na V\cdot\na(n-\bar{n})dx\bigg|
  \le \delta\int_\Omega\big(|\na n|^2+|\na p|^2+|\na D|^2\big)dx
  + C(\delta,L_*),
\end{align*}
where $\delta>0$ is arbitrary. Similar estimates hold for $p$ and $D$. Adding these inequalities and setting $Y:=\|n-\bar{n}\|_{L^2(\Omega)}^2 + \|p-\bar{p}\|_{L^2(\Omega)}^2 + \|D\|_{L^2(\Omega)}^2$, we obtain
\begin{align*}%\label{4.dYdt}
  \frac{dY}{dt}
  \le -\int_\Omega\big(|\na n|^2+|\na p|^2+|\na D|^2\big)dx + C(L_*)
  \le -cY + C\|D_I\|_{L^1(\Omega)}^2 + C(L_*),
\end{align*}
where we used the Poincar\'e inequality (for $n$ and $p$) and the Poincar\'e--Wirtinger inequality (for $D$). By Gronwall's lemma, $Y(t)\le \max\{Y(0),C(L_*)/c\}$, finishing the proof.
\end{proof}

Next, we derive a time-uniform bound for the Burg entropy. 

\begin{lemma}[Time-uniform Burg entropy]\label{lem.unifburg}
Let $(n,p,D,V)$ be a weak solution to \eqref{1.n}--\eqref{1.bcD} on $(0,T)$. Then 
\begin{align*}
  \int_\Omega(-\log D(t))dx 
  \le C(M,\|\na V\|_{L^\infty(0,T;L^2(\Omega))}),
  \quad M\ge\sup_{0<t<T}\|D(t)\|_{L^\infty(\Omega)}.
\end{align*}
\end{lemma}

\begin{proof}
If $\|\na V(t)\|_{L^2(\Omega)}\le L$ for all $t>0$, inequality \eqref{3.burg} yields $\int_\Omega(-\log D(t))dx\le \int_\Omega(-\log D_I)dx + \frac12 L^2t$, and the right-hand side is not bounded uniformly in time. To prevent the linear growth, we use the dissipation term $\frac12\|\na\log D\|_{L^2(\Omega)}^2$ instead of dropping it. To close the estimate
\begin{align}\label{5.burg}
  \frac{d}{dt}\int_\Omega(-\log D)dx 
  \le -\frac12\|\na\log D\|_{L^2(\Omega)}^2 + \frac{L^2}{2},
\end{align}
we need a Poincar\'e inequality for $\log D$, although $D$ has no boundary data. Mass conservation and the upper bound for $D$ imply that $D$ is not small on a set of positive measure, and the Poincar\'e inequality for functions vanishing on such sets applies. More precisely, let $\eps=\|D_I\|_{L^1(\Omega)}/(2|\Omega|)$ and $Z=\{D\ge\eps\}$. It follows from $\|D_I\|_{L^1(\Omega)}=\int_\Omega Ddx \le M|Z| + \eps|\Omega|$ that $|Z|\ge \|D_I\|_{L^1(\Omega)}/(2M)$, so $Z$ is a set of positive measure. The function $g=(\log\eps-\log D)_+$ vanishes on $Z$ and satisfies $|\na g|\le|\na\log D|$, thus $\|g\|_{L^2(\Omega)}\le C(M)\|\na\log D\|_{L^2(\Omega)}$ by the Poincar\'e inequality. The Burg entropy satisfies
\begin{align*}
  \int_\Omega(-\log D)dx \le \int_\Omega(g-\log\eps)dx
  \le C(M)\|\na\log D\|_{L^2(\Omega)} + C.
\end{align*}
Hence, estimating $\|\na\log D\|_{L^2(\Omega)}$ in \eqref{5.burg}, we find that
\begin{align*}
  \frac{d}{dt}\int_\Omega(-\log D)dx
  \le -c(M)\bigg(\int_\Omega(-\log D)dx - C\bigg)_+^2 + \frac{L^2}{2}.
\end{align*}
By the nonlinear Gronwall lemma, $\int_\Omega(-\log D)dx\le \max\{\int_\Omega(-\log D_I)dx,C(M,L)\}$, where the first argument satisfies $\int_\Omega(-\log D_I)dx\le |\Omega||\log m_I|$, finishing the proof.
\end{proof}

\begin{proposition}[Bounds depending on $L\log L$ norms]
\label{prop.unifbd}
Let $(n,p,D,V)$ be a weak solution to \eqref{1.n}--\eqref{1.bcD} on $(0,T)$ and assume that $L(t)\le L_*$ for $0<t<T$, recalling definition \eqref{5.L} of $L(t)$. Then there exist constants $m(L_*),M(L_*)>0$, depending only on $L_*$ and the data, but not on $T$, such that
\begin{align*}
  m(L_*)\le n,p,D\le M(L_*)\quad\mbox{in }\Omega_T, \quad
  \sup_{0<t<T}\|V(t)\|_{W^{1,r}(\Omega)}\le C(L_*).
\end{align*}
\end{proposition}

\begin{proof}
We only sketch the proof, since it is similar to the results in Section \ref{sec.bd}. The $L\log L$ bound gives an $L^2(\Omega)$ bound (Lemma \ref{lem.L2H1}). By elliptic regularity, the $L^2(\Omega)$ bound gives $\na V\in L^r(\Omega)$ with $r>2$ uniformly in time. The upper bounds follow from an Alikakos iteration as in Lemma \ref{lem.up}, but now in differential form and with the weight $e^t$, which makes the constants independent of $T$. 

The proof of time-uniform lower bounds for $n$ and $p$ differs from that of Lemma \ref{lem.lownp}, as the latter yields only the lower bound $m_0e^{-\mu t}$, which converges to zero as $t\to\infty$. To obtain a time-uniform lower bound, the source term must dominate the drift term on $\{n<m\}$, which requires that $\min_{[0,M]^2}R_0>0$ from Assumption (A3). This condition is used only at this point. Indeed, arguing as in the proof of Lemma \ref{lem.lownp}, since $p\le M$ and $(n-m)_-\le 0$,
\begin{align*}
  \frac12&\frac{d}{dt}\int_\Omega(n-m)_-^2 dx
  + \int_\Omega|\na(n-m)_-|^2 dx
  \le -Cm\int_\Omega(n-m)_- dx \\
  &\phantom{xx}+ \int_\Omega R_0(n,p)(n-m)_- dx 
  - \int_\Omega R_0(n,p)np(n-m)_- dx \\
  &\le \Big(Cm-\min_{[0,M]^2}R_0 - MmC\Big)\int_\Omega|(n-m)_-|dx \le C,
\end{align*}
if we choose $m>0$ sufficiently small. The lower bound for $D$ follows from the time-uniform Burg entropy bound and an Alikakos iteration for $-(\log D+m)_-$, again with the weight $e^t$, similarly as in the proofs of Lemmas \ref{lem.burg} and \ref{lem.lowD}.
\end{proof}

The proof of Theorem \ref{thm.bd} is now a consequence of Proposition \ref{prop.feb}, giving $\sup_{t>0}H(t)<\infty$, Lemma \ref{lem.cons}, showing that $\sup_{t>0}L(t)<\infty$, and Proposition \ref{prop.unifbd}. 

%%%%%%%%%%%%%%%%%%%%%%%%%%%%%%%%%%%%%%%%%%%%%%%%%%%%%%%%%%%%

\section{Exponential decay to thermal equilibrium}\label{sec.exp1}

We show that the solution $(n,p,D,V)$ to \eqref{1.n}--\eqref{1.bcD} converges exponentially fast to the thermal equilibrium $(n_{\rm eq},p_{\rm eq},D_{\rm eq},V_{\rm eq})$ under the assumptions of Theorem \ref{thm.exp1}. 

\subsection{Thermal equilibrium}%\label{sec.eq}

The thermal equilibrium is a solution to the integro--differen\-tial problem \eqref{1.eq}--\eqref{1.cD}. We show that such a solution exists. 

\begin{proposition}[Existence of thermal equilibrium]\label{prop.Veq}
There exists a unique weak solution $V_{\rm eq}\in H^1(\Omega)\cap L^\infty(\Omega)$ to \eqref{1.Veq}--\eqref{1.cD}. 
\end{proposition}

\begin{proof}
The result is proved by minimization of the functional (as in, e.g., \cite[Lemma 2.5]{ABCHM25})
\begin{align*}
  J(u) = \frac{\lambda^2}{2}\int_\Omega|\na u|^2 dx
  + \int_\Omega(c_n e^u + c_p e^{-u})dx
  + \|D_I\|_{L^1(\Omega)}\log\bigg(\int_\Omega e^{-u}dx\bigg)
  + \int_\Omega Audx
\end{align*}
for $u\in H^1(\Omega)$, with $J(u)=+\infty$ if $e^u\not\in L^1(\Omega)$ or $e^{-u}\not\in L^1(\Omega)$. Its Euler--Lagrange equation on $U:=H_D^1(\Omega)+\bar{V}_{\rm eq}$ is the weak formulation of \eqref{1.Veq}. The logarithm of the integral accounts for the mass constraint. To minimize $J$ on $U$, we first prove some properties.

We claim that $J$ is convex on $H^1(\Omega)$. Since $u\mapsto\int_\Omega(c_n e^u+c_p e^{-u})dx$ is strictly convex, it remains to verify the convexity of $F(u)=\log(\int_\Omega e^{-u}dx)$. Indeed, by H\"older's inequality with exponents $1/\theta$ and $1/(1-\theta)$,
\begin{align*}
  F(\theta u_1&+(1-\theta)u_2)
  = \log\bigg(\int_\Omega e^{-(\theta u_1+(1-\theta)u_2)}dx\bigg) \\
  &\le \log\bigg\{\bigg(\int_\Omega e^{-u_1}dx\bigg)^\theta
  \bigg(\int_\Omega e^{-u_2}dx\bigg)^{1-\theta}\bigg\}
  = \theta F(u_1)+(1-\theta)F(u_2).
\end{align*}
The sum of a strictly convex and a convex function is strictly convex, so $J$ is strictly convex on $\{J<\infty\}$. The weak lower semicontinuity of $J$ is a consequence of the weak lower semicontinuity of the Dirichlet integral and the lower semicontinuity of the exponential term. We prove the coercivity of $J$. Let $c=\min\{c_n,c_p\}$. Then $c_n e^u+c_p e^{-u}\ge c e^{|u|}\ge cu^2/2$. We deduce from Jensen's inequality applied to the exponential that
\begin{align}\label{4.jensen}
  \frac{1}{|\Omega|}\int_\Omega e^{-u}dx 
  \ge \exp\bigg({-}\frac{1}{|\Omega|}\int_\Omega udx\bigg).
\end{align}
Hence, $F(u)\ge \log|\Omega|-|\Omega|^{-1}\int_\Omega udx\ge\log|\Omega|-|\Omega|^{-1/2}\|u\|_{L^2(\Omega)}$. This leads to
\begin{align*}
  J(u) \ge \frac{\lambda^2}{2}\|\na u\|_{L^2(\Omega)}^2
  + \frac{c}{4}\|u\|_{L^2(\Omega)}^2 - C(A,|\Omega|).
\end{align*}

Since $J(\bar{V}_{\rm eq})<\infty$, we can apply the direct method of the calculus of variations to conclude a minimizer $V_{\rm eq}$ of $J$ on $U$, which is unique by the strict convexity of $J$. Then the first variation of $J$ at $V_{\rm eq}$ vanishes, which yields the weak formulation of \eqref{1.Veq}. 
\end{proof}

We prove $L^\infty(\Omega)$ bounds for $V_{\rm eq}$ and $(n_{\rm eq},p_{\rm eq},D_{\rm eq})$.

\begin{lemma}[Bounds for $V_{\rm eq}$]
Let $V_{\rm eq}\in H_D^1(\Omega)+\bar{V}_{\rm eq}$ be the unique minimizer of $J$ from Proposition \ref{prop.Veq}. Then there exist constants 
\begin{itemize}
\item[\rm (i)] $C_1>0$ depending on $\|D_I\|_{L^1(\Omega)}$ and $\Omega$ such that $0<c_D[V_{\rm eq}]\le C_1$;
\item[\rm (ii)] $C_2>0$ depending on the data such that $\|V_{\rm eq}\|_{L^\infty(\Omega)}\le C_2$;
\item[\rm (iii)] $m_{\rm eq},M_{\rm eq}>0$ depending on $C_1$, $C_2$ such that
\begin{align*}
  m_{\rm eq}\le n_{\rm eq},\,p_{\rm eq},\,D_{\rm eq}\le M_{\rm eq}.
\end{align*}
\end{itemize}
\end{lemma}

\begin{proof}
(i) The coercivity of $J$ immediately implies an $H^1(\Omega)$ bound for $V_{\rm eq}$ depending on $J(\bar{V}_{\rm eq})$. As in \eqref{4.jensen}, we have
\begin{align*}
  \int_\Omega e^{-V_{\rm eq}}dx
  \ge |\Omega|\exp\bigg(-\frac{1}{|\Omega|}\int_\Omega V_{\rm eq}
  dx\bigg)\ge |\Omega|
  \exp\bigg(-\frac{\|V_{\rm eq}\|_{L^2(\Omega)}}{|\Omega|^{1/2}}\bigg),
\end{align*}
which yields an upper bound for $c_D[V_{\rm eq}]$. (ii) Let $K>\|\bar{V}_{\rm eq}\|_{L^\infty(\Omega)}$ and $m\in\N$. We use the test function $\phi_m=\min\{m,(V_{\rm eq}-K)_+\}$ in the weak formulation of \eqref{1.Veq}:
\begin{align*}
  \lambda^2&\int_{\{K<V_{\rm eq}<m+K\}}|\na V_{\rm eq}|^2 dx
  = -\int_{\{V_{\rm eq}>K\}}\big(c_n e^{V_{\rm eq}} 
  - c_p e^{-V_{\rm eq}} - c_D[V_{\rm eq}]e^{-V_{\rm eq}} 
  + A\big)\phi_m dx \\
  &\le -\int_{\{V_{\rm eq}>K\}}\big(c_n e^K - c_p e^{-K}
  - C_1 - \|A\|_{L^\infty(\Omega)}\big)\phi_m dx \le 0,
\end{align*}
if $K$ is sufficiently large. Monotone convergence implies for $m\to\infty$ that $\int_\Omega|\na V_{\rm eq}|^2$ $1_{\{V_{\rm eq}>K\}}dx=0$ and hence $\na(V_{\rm eq}-K)_+=0$ or $V_{\rm eq}\le K$ in $\Omega$. The lower bound is shown in a similar way using the test function $\max\{-m,(V_{\rm eq}+K)_-\}$. (iii) The bounds follow from definition \eqref{1.eq} of $n_{\rm eq}$, $p_{\rm eq}$, and $D_{\rm eq}$. 
\end{proof}

%%%%%%%%%%%%%%%%%%%%

\subsection{Harmonic extensions}

Since the equilibrium densities may not satisfy the Dirichlet conditions, we need to introduce boundary-adapted equilibrium densities. To this end, let $E:H^1(\Omega)\to H^1(\Omega)$ be the harmonic extension operator, i.e., $w=E[u]$ is the unique weak solution to $\Delta w=0$ in $\Omega$, $w=u$ on $\Gamma_D$, and $\na w\cdot\nu=0$ on $\Gamma_N$. It satisfies $\|E[u]\|_{H^1(\Omega)}\le C\|u\|_{H^1(\Omega)}$, $\|E[u]\|_{L^\infty(\Omega)}\le\|u\|_{L^\infty(\Omega)}$, and $\|\pa_t(E[u])\|_{L^\infty(\Omega)}\le \|\pa_t u\|_{L^\infty(\Omega)}$ for suitable functions $u$. The boundary-adapted equilibrium densities are defined by
\begin{align*}
  \widehat{n} = n_{\rm eq}e^N, \quad \widehat{p} = p_{\rm eq}e^P,
\end{align*}
where $N = E[\log(\bar{n}/n_{\rm eq})]$, $P = E[\log(\bar{p}/p_{\rm eq})]$ are bounded in $L^\infty(\Omega)$. We also introduce the potential $U = E[\bar{V}-\bar{V}_{\rm eq}]$, which equals $\bar{V}-\bar{V}_{\rm eq}$, since both functions are harmonic (Assumption (A5)). The function $\bar V_{\rm eq}$ is (the trace of) a harmonic function with $\na\bar V_{\rm eq}\cdot\nu=0$ on $\Gamma_N$ in the weak sense, e.g.\ the limit of $\bar V(t)$. We also define the relative free energy 
\begin{align*}
  H_{\rm rel} = \int_\Omega\big(h(n|\widehat{n}) + h(p|\widehat{p})
  + h(D|D_{\rm eq})\big)dx + \frac{\lambda^2}{2}\int_\Omega
  |\na W|^2 dx,
\end{align*}
where $W = V-V_{\rm eq}-U$ and $h(u|\widehat{u}) = u(\log(u/\widehat{u})-1) + \widehat{u}$. The idea of this construction is that $\widehat{n}=\bar{n}$, $\widehat{p}=\bar{p}$, $W=0$ (i.e., $V=V_{\rm eq}+U$) on $\Gamma_D$, so that the shifted quasi-Fermi potentials
\begin{align*}
  \phi_n = \log\frac{n}{\widehat{n}}-W, \quad
  \phi_p = \log\frac{p}{\widehat{p}}+W, \quad
  \phi_D= \log\frac{D}{D_{\rm eq}}+W
\end{align*}
satisfy $\phi_n=0$ and $\phi_p=0$ on $\Gamma_D$ and can be used as test functions.

We wish to measure the distance of the boundary data from the equilibrium data by
\begin{align}\label{6.S}
  \mathcal S(t) = \|N\|_{H^1(\Omega)}^2 + \|P\|_{H^1(\Omega)}^2
  + \|\pa_t N\|_{L^2(\Omega)}^2 + \|\pa_t P\|_{L^2(\Omega)}^2
  + \|U\|_{H^1(\Omega)}^2.
\end{align}

\begin{lemma}%\label{lem.S}
There exists $C>0$ such that for a.e.\ $t>0$,
\begin{align*}
  \mathcal S(t) \le C\bigg(\Lambda(t) + \bar\Lambda(t)
  + \bigg(\int_t^\infty\bar\Lambda(s)ds\bigg)^2\bigg), \quad
  \bar\Lambda = \Lambda + \|\pa_t\bar{V}\|_{H^1(\Omega)},
\end{align*}
recalling definition \eqref{1.Lambda} of $\Lambda$. 
\end{lemma}

\begin{proof}
We only give a sketch of the proof, since we may directly assume a decay property for $\mathcal{S}(t)$. Recalling that $\bar{\mu}_n = \log\bar{n}-\bar{V}$, $\bar{\mu}_p = \log\bar{p}+\bar{V}$ from Assumption (A5) and $\log n_{\rm eq} = \log c_n +  V_{\rm eq}$, $\log p_{\rm eq} = \log c_p -  V_{\rm eq}$ by definition of the thermal equilibrium, it follows that 
\begin{align*}
  N = E[\log\bar{n} - \log n_{\rm eq}]
  = E[\bar{\mu}_n - \log c_n + U + (\bar{V}_{\rm eq}-V_{\rm eq})]
  = E[\bar{\mu}_n-\log c_n] + U
\end{align*}
and similarly, $P = E[\bar{\mu}_p-\log c_p] - U$. We infer that
\begin{align*}
  \|N\|_{H^1(\Omega)} \le C\|\bar{\mu}_n-\log c_n\|_{H^1(\Omega)}
  + \|U\|_{H^1(\Omega)}, \quad 
  \|U\|_{H^1(\Omega)} \le C\|\bar{V}-\bar{V}_{\rm eq}\|_{H^1(\Omega)},
\end{align*}
and similar for the $H^1(\Omega)$ norm of $P$. Together with $\|\pa_t N\|_{L^2(\Omega)} \le C\|\pa_t\bar{\mu}_n\|_{L^\infty(\Omega)} + \|\pa_t\bar{V}\|_{L^2(\Omega)}$, this gives
\begin{align*}
  \mathcal S(t) &\le C\|\bar{\mu}_n-\log c_n\|_{H^1(\Omega)}^2
  + C\|\bar{\mu}_p-\log c_p\|_{H^1(\Omega)}^2
  + C\|\pa_t\bar{\mu}_n\|_{L^\infty(\Omega)}^2 \\
  &\phantom{xx}+ C\|\pa_t\bar{\mu}_p\|_{L^\infty(\Omega)}^2
  + C\|\bar{V}-\bar{V}_{\rm eq}\|_{H^1(\Omega)}^2
  + C\|\pa_t\bar{V}\|_{L^2(\Omega)}^2.
\end{align*}
It can be shown that the right-hand side is estimated according to
\begin{align*}
  &\|\bar{V}(t)-\bar{V}_{\rm eq}\|_{H^1(\Omega)}^2
  \le \bigg(\int_t^\infty\bar\Lambda(\tau)d\tau\bigg)^2, \quad
  \|\na\bar{\mu}_n\|_{L^\infty(\Omega)}^2\le C\Lambda, 
  \quad \|\pa_t\bar{\mu}_n\|_{L^\infty(\Omega)}\le C\bar\Lambda, \\
  & \|\bar{\mu}_n-\log c_n\|_{H^1(\Omega)}^2
  \le C\big(\|\bar{\mu}_n-\log c_n\|_{L^\infty(\Omega)}^2
  + \|\na\bar{\mu}_n\|_{L^\infty(\Omega)}^2\big) \\
  &\phantom{\|\bar{\mu}_n-\log c_n\|_{H^1(\Omega)}^2}
  \le C\bigg\{\bigg(\int_t^\infty\Lambda(s)ds\bigg)^2 
  + \Lambda(t)\bigg\},
\end{align*}
and the same bounds hold for $\na\bar{\mu}_p$, $\pa_t\bar{\mu}_n$, and $\bar{\mu}_p-\log c_p$. Observe that these estimates use that $\log c_n=\lim_{t\to\infty}\bar\mu_n(t)$ and $\bar V_{\rm eq}=\lim_{t\to\infty}\bar V(t)$, e.g.\ $\|\bar\mu_n(t)-\log c_n\|_{L^\infty}\le\int_t^\infty\|\pa_s\bar\mu_n\|_{L^\infty}ds$. The proof is complete.
\end{proof}

%%%%%%%%%%%%%%%%%%%%

\subsection{Exponential decay to the thermal equilibrium}

We split the proof of Theorem \ref{thm.exp1} into several steps.

{\em Step 1: Relative energy inequality.} We compute
\begin{align*}
  \frac{dH_{\rm rel}}{dt} = \langle\pa_t n,\phi_n\rangle
  + \langle\pa_t p,\phi_p\rangle + \langle\pa_t D,\phi_D\rangle
  - \int_\Omega \big((n-\widehat{n})\pa_tN 
  + (p-\widehat{p})\pa_tP\big)dx.
\end{align*}
Taking into account that $\na(\log n-V) = \na(\phi_n + (\log n_{\rm eq}-V_{\rm eq}) + (N-U)) = \na\phi_n + \na(N-U)$, Young's inequality, and the $L^\infty(\Omega)$ bound from Proposition \ref{prop.unifbd} give
\begin{align*}
  \langle\pa_t n,\phi_n\rangle 
  &= -\int_\Omega n\na(\log n-V)\cdot\na\phi_n dx
  - \int_\Omega R(n,p)\phi_n dx \\
  &= -\int_\Omega n|\na\phi_n|^2 dx 
  - \int_\Omega n\na(N-U)\cdot\na\phi_n dx
  - \int_\Omega R(n,p)\phi_n dx \\
  &\le -\frac12\int_\Omega n|\na\phi_n|^2 dx 
  + \frac{M}{2}\int_\Omega|\na(N-U)|^2 dx
  - \int_\Omega R(n,p)\phi_n dx,
\end{align*}
where $M$ is an upper bound for $n$, and similarly for $p$ and $D$. Using $n_{\rm eq}p_{\rm eq}=1$ and consequently $\log(\widehat{n}\widehat{p})=N+P$, the generation--recombination terms combine to
\begin{align*}
  -\int_\Omega& R(n,p)(\phi_n+\phi_p)dx
  = -\int_\Omega R_0(n,p)(np-1)\big(\log(np)-(N+P)\big)dx \\
  &= \int_\Omega R_0(n,p)(1-np)\log(np)dx 
  + \int_\Omega R_0(n,p)(np-1)(N+P)dx.
\end{align*}
A Taylor expansion leads to $(1-z)\log z\le -c(K)(1-z)^2$ for all $0<z\le K$ and some $c(K)>0$. Hence, with $z=np$, this shows that
\begin{align*}
  -\int_\Omega R_0(n,p)(np-1)\log(np)dx 
  \le -c(M^2)\int_\Omega R_0(n,p)(1-np)^2 dx,
\end{align*}
and together with $(np-1)(N+P)\le \delta(1-np)^2 + C(\delta)(N+P)^2$ for $\delta>0$ and the boundedness of $R_0$,
\begin{align*}
  -\int_\Omega R(n,p)(\phi_n+\phi_p)dx
  \le -\frac12c(M^2)\int_\Omega R_0(n,p)(1-np)^2 dx
  + C\int_\Omega(N+P)^2 dx.
\end{align*}
Summarizing, we obtain
\begin{align}\label{4.dHdt}
  & \frac{dH_{\rm rel}}{dt} + \frac12\mathcal{D}
  + c\int_\Omega R_0(n,p)(1-np)^2 dx
  \le C\mathcal{S}(t) \\
  &\phantom{dH_{\rm reli}}+ \|\pa_t N\|_{L^2(\Omega)}
  \|n-\widehat{n}\|_{L^2(\Omega)}
  + \|\pa_t P\|_{L^2(\Omega)}\|p-\widehat{p}\|_{L^2(\Omega)},
  \nonumber \\
  &\mbox{where}\quad\mathcal{D} = \int_\Omega\big(n|\na\phi_n|^2
  + p|\na\phi_p|^2 + D|\na\phi_D|^2\big)dx. \nonumber 
\end{align}

{\em Step 2: Energy dissipation inequality.} It follows from $\lambda^2\Delta W = \lambda^2\Delta(V-V_{\rm eq}) = (n-n_{\rm eq}) - (p-p_{\rm eq}) - (D-D_{\rm eq})$ that
\begin{align*}
  I &:= \int_\Omega\big((n-\widehat{n})\phi_n
  + (p-\widehat{p})\phi_p + (D-D_{\rm eq})\phi_D\big)dx \\
  &= \int_\Omega\bigg((n-\widehat{n})\log\frac{n}{\widehat{n}}
  + (p-\widehat{p})\log\frac{p}{\widehat{p}}
  + (D-D_{\rm eq})\log\frac{D}{D_{\rm eq}}\bigg)dx \\
  &\phantom{xx}- \int_\Omega\big((n-p-D)
  - (\widehat{n}-\widehat{p}-D_{\rm eq})\big) W dx \\
  &= \int_\Omega\bigg((n-\widehat{n})\log\frac{n}{\widehat{n}}
  + (p-\widehat{p})\log\frac{p}{\widehat{p}}
  + (D-D_{\rm eq})\log\frac{D}{D_{\rm eq}}\bigg)dx \\
  &\phantom{xx}+ \lambda^2\int_\Omega|\na W|^2 dx
  + \int_\Omega\big((\widehat{n}-n_{\rm eq}) 
  - (\widehat{p}-p_{\rm eq})\big)Wdx.
\end{align*}
The aim is to estimate this expression from below and above. On the one hand, since $(a-b)\log(a/b)\ge (a-b)^2/\max\{a,b\}$ and $|\widehat{n}-n_{\rm eq}| = n_{\rm eq}|e^N-1|\le C|N|$ (and similarly for $p$), the Poincar\'e inequality for $W\in H_D^1(\Omega)$ yields
\begin{align*}
  I &\ge c\big(\|n-\widehat{n}\|_{L^2(\Omega)}^2 
  + \|p-\widehat{p}\|_{L^2(\Omega)}^2 
  + \|D-D_{\rm eq}\|_{L^2(\Omega)}^2\big)
  + \frac{\lambda^2}{2}\|\na W\|_{L^2(\Omega)}^2 \\
  &\phantom{xx}- C\int_\Omega(|N|^2+|P|^2)dx. 
\end{align*}
On the other hand, we have $\phi_n,\phi_p\in H_D^1(\Omega)$, and $\int_\Omega(D-D_{\rm eq})dx=0$ allows us to subtract the mean value of $\phi_D$. Then we can apply the Poincar\'e and Poincar\'e--Wirtinger inequalities, and together with the time-uniform lower bound $n,p,D\ge m$ from Proposition \ref{prop.unifbd}, 
\begin{align*}
  I &\le \frac{c}{2}\big(\|n-\widehat{n}\|_{L^2(\Omega)}^2 
  + \|p-\widehat{p}\|_{L^2(\Omega)}^2 
  + \|D-D_{\rm eq}\|_{L^2(\Omega)}^2\big) \\
  &\phantom{xx}+ C\big(\|\na\phi_n\|_{L^2(\Omega)}^2
  + \|\na\phi_p\|_{L^2(\Omega)}^2 
  + \|\na\phi_D\|_{L^2(\Omega)}^2\big) \\
  &\le \frac{c}{2}\big(\|n-\widehat{n}\|_{L^2(\Omega)}^2 
  + \|p-\widehat{p}\|_{L^2(\Omega)}^2 
  + \|D-D_{\rm eq}\|_{L^2(\Omega)}^2\big) 
  + \frac{C}{m}\mathcal{D}.
\end{align*}
We combine the lower and upper bounds for $I$:
\begin{align}\label{4.aux}
  \frac{c}{2}\big(\|n-\widehat{n}\|_{L^2(\Omega)}^2 
  &+ \|p-\widehat{p}\|_{L^2(\Omega)}^2 
  + \|D-D_{\rm eq}\|_{L^2(\Omega)}^2\big)
  + \frac{\lambda^2}{2}\|\na W\|_{L^2(\Omega)}^2 \\
  &\le \frac{C}{m}\mathcal{D}
  + C\big(\|N\|_{L^2(\Omega)}^2 + \|P\|_{L^2(\Omega)}^2\big).
  \nonumber 
\end{align}
Furthermore, we deduce from $h(a|b)\le (a-b)^2/(2\min\{a,b\})$ that
\begin{align} \label{4.H}
  H_{\rm rel} &\le C\big(\|n-\widehat{n}\|_{L^2(\Omega)}^2 
  + \|p-\widehat{p}\|_{L^2(\Omega)}^2 
  + \|D-D_{\rm eq}\|_{L^2(\Omega)}^2\big)
  + \frac{\lambda^2}{2}\|\na W\|_{L^2(\Omega)}^2 \\
  &\le C\mathcal{D} + C\big(\|N\|_{L^2(\Omega)}^2 
  + \|P\|_{L^2(\Omega)}^2\big). \nonumber 
\end{align}
By Young's inequality and \eqref{4.aux}, the last two terms in \eqref{4.dHdt} are bounded by $\frac14\mathcal{D} + C\|\pa_t N\|_{L^2(\Omega)}^2$ $+C\|\pa_t P\|_{L^2(\Omega)}^2$. We infer from \eqref{4.dHdt}, also using \eqref{4.H}, that 
\begin{align*}
  & \frac{dH_{\rm rel}}{dt} \le -\frac14\mathcal{D}
  + C\mathcal{S} \le -c H_{\rm rel} + C\mathcal{S},
\end{align*}
recalling definition \eqref{6.S} of $\mathcal S$. 

{\em Step 3: End of the proof.} 
We deduce from the decay assumption \eqref{1.A6} that $\frac{d}{dt}H_{\rm rel} + cH_{\rm rel}\le Ce^{-\kappa_0 t}$. By Gronwall's lemma, 
\begin{align*}
  H_{\rm rel}(t) \le e^{-ct}H_{\rm rel}(0)
  + C\int_0^t e^{-c(t-s)-\kappa_0 s}ds.
\end{align*}
We estimate the last term:
\begin{align*}
  \int_0^t e^{-c(t-s)-\kappa_0 s}ds
  &= \begin{cases}
  te^{-ct} &\mbox{if }c=\kappa_0 \\
  \dfrac{e^{-\kappa_0 t}-e^{-ct}}{c-\kappa_0} &\mbox{if }c\neq\kappa_0
  \end{cases} \\
  &\le \min\bigg\{t,\frac{1}{|c-\kappa_0|}\bigg\}
  e^{-\min\{c,\kappa_0\}t} \le C(\kappa)e^{-\kappa t}
\end{align*}
for any $0<\kappa<\min\{c,\kappa_0\}$. 

Finally, we prove the $L^2(\Omega)$ decay bounds. It follows from $h(a|b)\ge(a-b)^2/(2\max\{a,b\})$ that
\begin{align*}
  \|n-\widehat{n}\|_{L^2(\Omega)}^2 
  + \|p-\widehat{p}\|_{L^2(\Omega)}^2
  + \|D-D_{\rm eq}\|_{L^2(\Omega)}^2 \le 2MH_{\rm rel}[n,p,D,V]. 
\end{align*}
We combine this with 
\begin{align*}
  \|\widehat{n}-n_{\rm eq}\|_{L^2(\Omega)}
  = \|n_{\rm eq}(1-e^N)\|_{L^2(\Omega)}
  \le \|n_{\rm eq}\|_{L^\infty(\Omega)}
  e^{\|N\|_{L^\infty(\Omega)}}\|N\|_{L^2(\Omega)}
  \le Ce^{-\kappa_0 t/2},
\end{align*}
and similarly for $p$. Then
\begin{align*}
  \|n-n_{\rm eq}\|_{L^2(\Omega)}
  + \|p-p_{\rm eq}\|_{L^2(\Omega)} + \|D-D_{\rm eq}\|_{L^2(\Omega)}
  \le C(H_{\rm rel}^{1/2} + e^{-\kappa_0 t/2}) \le Ce^{-\kappa t/2}.
\end{align*}
For the potential, we have
\begin{align*}
  \|V-V_{\rm eq}\|_{H^1(\Omega)}
  = \|W+U\|_{H^1(\Omega)} \le C\|\na W\|_{L^2(\Omega)}
  + \|U\|_{H^1(\Omega)} \le Ce^{-\kappa t/2}.
\end{align*}
This completes the proof of Theorem \ref{thm.exp1}.

%%%%%%%%%%%%%%%%%%%%%%%%%%%%%%%%%%%%%%%%%%%%%%%%%%%%%%%%%%%%%

\section{Exponential decay to steady states with non-vanishing fluxes}
\label{sec.exp2}

We recall that, by definition, $(n^*,p^*,D^*,V^*)\in H^1(\Omega)^4$ is a steady state if the stationary versions of \eqref{1.n}--\eqref{1.V}, the boundary conditions \eqref{1.Dir}--\eqref{1.bcD}, and $m^*\le n^*,p^*,D^*\le M^*$ hold in $\Omega$ for some $0<m^*\le M^*$. Since the stationary charge $n^*-p^*-D^*+A$ is bounded, we have
\begin{align*}
  \|V^*\|_{W^{1,r}(\Omega)}\le C,
\end{align*}
where $C>0$ only depends on $m^*$, $M^*$, and the data. Then $\Delta n^*=\diver(n^*\na V^*)+R(n^*,p^*)\in W_D^{-1,r}(\Omega)$, and the regularity result in \cite{Gro89} implies that $n^*,p^*\in W^{1,r}(\Omega)$ with $r>2$. Thus, the non-equilibrium parameter
\begin{align}\label{5.eta}
  \eta := \|\na(\log n^*-V^*)\|_{L^r(\Omega)}
  + \|\na(\log p^*+V^*)\|_{L^r(\Omega)}
  + \|1-n^*p^*\|_{L^\infty(\Omega)}
\end{align}
is well defined. The stationary electron and hole current densities $J_n^*=n^*\na(\log n^*-V^*)$ and $J_p^*=p^*\na(\log p^*+V^*)$ are of the order $\eta$ but do not vanish in general. The thermal equilibrium of the previous section is a steady state with $\eta=0$. We expect that the uniqueness of the steady state for sufficiently small $\eta$ can be proved by the implicit function theorem similarly as in \cite{NaWo91}.

Interestingly, stationary vacancies are already in equilibrium. Indeed, the test function $\log D^*+V^*$ in the weak formulation of the stationary version of \eqref{1.D} yields $\int_\Omega D^*|\na(\log D^*+V^*)|^2 dx = 0$, and because of the strict positivity of $D^*$, $\log D^*+V^*=\mbox{const.}$, where the constant is determined by the mass, i.e.
\begin{align*}
  D^* = \|D_I\|_{L^1(\Omega)}\bigg(\int_\Omega e^{-V^*}dx\bigg)^{-1}
  e^{-V^*}.
\end{align*}
Thus, the steady-state problem is the stationary bipolar model with an additional nonlocal term in the Poisson equation, similarly as in \eqref{1.Veq}. This term is the derivative of the convex functional $\log\int_\Omega e^{-V^*}dx$, up to a factor, and hence monotone. This may allow for the proof of a unique weak solution if $\eta>0$ is sufficiently small.

Let $(n,p,D,V)$ be the weak solution to \eqref{1.n}--\eqref{1.bcD} with the same boundary data and the same total vacancy mass as the steady state. We define the relative entropy with respect to the steady state and the corresponding dissipation:
\begin{align*}
  H^* &= \int_\Omega\big(h(n|n^*) + h(p|p^*) + h(D|D^*)\big)dx
  + \frac{\lambda^2}{2}\int_\Omega|\na(V-V^*)|^2 dx, \\
  \mathcal{D^*} &= \int_\Omega\big(n|\na(\mu_n-\mu_n^*)|^2
  + p|\na(\mu_p-\mu_p^*)|^2 + D|\na(\mu_D-\mu_D^*)|^2\big)dx,
\end{align*}
where $\mu_n=\log n-V$, $\mu_p=\log p+V$, $\mu_D=\log D+V$, similarly for $\mu_n^*$, $\mu_p^*$, $\mu_D^*$, and recalling that $h(u|\widehat{u}) = u(\log(u/\widehat{u})-1) + \widehat{u}$. Notice that $\mu_n-\mu_n^*=\mu_p-\mu_p^*=0$ on $\Gamma_D$, since the boundary values of $(n,p,V)$ and $(n^*,p^*,V^*)$ coincide by assumption. 

\begin{lemma}[Relative free energy identity]\label{lem.Hstar}
Let $(n,p,D,V)$ be a weak solution to \eqref{1.n}--\eqref{1.bcD}. Then
\begin{align}\label{5.dHstar}
  \frac{dH^*}{dt} + \mathcal{D}^*
  &= -\int_\Omega\big((n-n^*)\na\mu_n^*\cdot\na(\mu_n-\mu_n^*)
  + (p-p^*)\na\mu_p^*\cdot\na(\mu_p-\mu_p^*)\big)dx \\
  &\phantom{xx}- \int_\Omega(R(n,p)-R(n^*,p^*))\log\frac{np}{n^*p^*}dx.
  \nonumber 
\end{align}
\end{lemma}

\begin{proof}
We compute
\begin{align}\label{5.aux}
  \frac{dH^*}{dt} = \langle\pa_t n,\mu_n-\mu_n^*\rangle
  + \langle\pa_t p,\mu_p-\mu_p^*\rangle
  + \langle\pa_t D,\mu_D-\mu_D^*\rangle.
\end{align}
Since $\diver(n^*\na\mu_n^*)=R(n^*,p^*)$, we obtain
\begin{align*}
  \langle\pa_t n,\mu_n-\mu_n^*\rangle
  &= -\int_\Omega n|\na(\mu_n-\mu_n^*)|^2 dx
  - \int_\Omega(n-n^*)\na\mu_n^*\cdot\na(\mu_n-\mu_n^*)dx \\
  &\phantom{xx}- \int_\Omega(R(n,p)-R(n^*,p^*))(\mu_n-\mu_n^*)dx,
\end{align*}
and similarly for the term involving $p$, while the term involving $D$ simplifies, because of $\na\mu_D^*=0$,
\begin{align*}
  \langle\pa_t D,\mu_D-\mu_D^*\rangle
 = -\int_\Omega D|\na(\mu_D-\mu_D^*)|^2 dx.
\end{align*}
The generation--recombination terms add up to 
\begin{align*}
  -\int_\Omega (R(n,p)-R(n^*,p^*))(\mu_n+\mu_p-\mu_n^*-\mu_p^*)dx
  = -\int_\Omega (R(n,p)-R(n^*,p^*))\log\frac{np}{n^*p^*}dx.
\end{align*}
Inserting these expressions into \eqref{5.aux} finishes the proof.
\end{proof}

\begin{lemma}\label{lem.rhs}
Let $(n,p,D,V)$ be a weak solution to \eqref{1.n}--\eqref{1.bcD} such that $0<m\le n,p,D\le M$ in $\Omega$, and let $\eta\le 1$. Then there exists a constant $C_1>0$ independent of $T$ such that $H^*\le C_1\mathcal{D}^*$ and 
\begin{align}\label{5.ineq}
  \frac{dH^*}{dt} + \mathcal{D}^*
  \le C_1\eta(H^* + \mathcal{D}^*).  
\end{align}
\end{lemma}

\begin{proof}
By Taylor expansion and the lower and upper bounds for the densities,
\begin{align*}
  c\big(&\|n-n^*\|_{L^2(\Omega)}^2 + \|p-p^*\|_{L^2(\Omega)}^2
  + \|D-D^*\|_{L^2(\Omega)}^2\big)
  + \lambda^2\|\na(V-V^*)\|_{L^2(\Omega)}^2 \\
  &\le \int_\Omega\bigg((n-n^*)\log\frac{n}{n^*}
  + (p-p^*)\log\frac{p}{p^*} + (D-D^*)\log\frac{D}{D^*}\bigg)dx \\
  &\phantom{xx}+ \lambda^2\|\na(V-V^*)\|_{L^2(\Omega)}^2 \\
  &= \int_\Omega\big((n-n^*)(\mu_n-\mu_n^*)
  + (p-p^*)(\mu_p-\mu_p^*) + (D-D^*)(\mu_D-\mu_D^*)\big)dx \\
  &\le \frac{c}{2}\big(\|n-n^*\|_{L^2(\Omega)}^2 
  + \|p-p^*\|_{L^2(\Omega)}^2 + \|D-D^*\|_{L^2(\Omega)}^2\big) \\
  &\phantom{xx}+ C\big(\|\na(\mu_n-\mu_n^*)\|_{L^2(\Omega)}^2
  + \|\na(\mu_p-\mu_p^*)\|_{L^2(\Omega)}^2
  + \|\na(\mu_D-\mu_D^*)\|_{L^2(\Omega)}^2\big),
\end{align*}
where the last step follows from the Young, Poincar\'e, and Poincar\'e--Wirtinger inequalities. Absorbing the first term on the right-hand side into the left-hand side, we find that
\begin{align*}
  \frac{c}{2}\big(\|n-n^*\|_{L^2(\Omega)}^2 
  + \|p-p^*\|_{L^2(\Omega)}^2 + \|D-D^*\|_{L^2(\Omega)}^2\big)
  + \lambda^2\|\na(V-V^*)\|_{L^2(\Omega)}^2
  \le Cm^{-1}\mathcal{D}^*.
\end{align*}
The inequality $H^*\le C_1\mathcal{D}^*$ follows from $h(n|n^*)\le C|n-n^*|^2$ and similarly for $p$ and $D$. 

For the proof of \eqref{5.ineq}, we first estimate $\|n-n^*\|_{H^1(\Omega)}$. By definition,
\begin{align}\label{5.aux2}
  \|\na(n-n^*)\|_{L^2(\Omega)}
  &\le \bigg\|n\na\log\frac{n}{n^*}\bigg\|_{L^2(\Omega)}
  + \|(n-n^*)\na\log n^*\|_{L^2(\Omega)} \\
  &= \|n\na(\mu_n-\mu_n^* + V - V^*)\|_{L^2(\Omega)}
  + \|(n-n^*)\na\log n^*\|_{L^2(\Omega)} \nonumber \\
  &\le M\|\na(\mu_n-\mu_n^*)\|_{L^2(\Omega)}
  + M\|\na(V-V^*)\|_{L^2(\Omega)} \nonumber \\
  &\phantom{xx} + (m^*)^{-1}\|\na n^*\|_{L^r(\Omega)}
  \|n-n^*\|_{L^{q}(\Omega)}, \nonumber 
\end{align} 
where $q=2r/(r-2)$. The regularity for the stationary problem gives $\|\na n^*\|_{L^r(\Omega)}\le C$. Furthermore, by the Gagliardo--Nirenberg inequality \eqref{1.GN} and the Young inequality,
\begin{align*}
  \|n-n^*\|_{L^q(\Omega)}
  \le \delta\|n-n^*\|_{H^1(\Omega)}
  + C(\delta)\|n-n^*\|_{L^2(\Omega)}.
\end{align*}
Absorbing the first term on the right-hand side into the left-hand side of \eqref{5.aux2}, 
\begin{align*}
  \|n-n^*\|_{H^1(\Omega)}
  &\le 2M\|\na(\mu_n-\mu_n^*)\|_{L^2(\Omega)}
  + 2M\|\na(V-V^*)\|_{L^2(\Omega)} \\
  &\phantom{xx}+ 2C\|n-n^*\|_{L^2(\Omega)}
  \le C\big((\mathcal{D}^*)^{1/2} + (H^*)^{1/2}\big).
\end{align*}
Then, using the H\"older and Gagliardo--Nirenberg inequalities as above, as well as definition \eqref{5.eta} of $\eta$,
\begin{align*}
  \bigg|\int_\Omega&(n-n^*)\na\mu_n^*\cdot\na(\mu_n-\mu_n^*)dx\bigg|
  \le \|n-n^*\|_{L^{q}(\Omega)}\|\na\mu_n^*\|_{L^r(\Omega)}
  \|\na(\mu_n-\mu_n^*)\|_{L^2(\Omega)} \\
  &\le C\eta m^{-1/2}\|n-n^*\|_{H^1(\Omega)}(\mathcal{D}^*)^{1/2}
  \le C(m)\eta\big((\mathcal{D}^*)^{1/2} + (H^*)^{1/2}\big)
  (\mathcal{D}^*)^{1/2} \\
  &\le C\eta(\mathcal{D}^* + H^*).
\end{align*}
Finally, the generation--recombination terms are written as
\begin{align*}
  -\int_\Omega&(R(n,p)-R(n^*,p^*))\log\frac{np}{n^*p^*}dx \\
  &= \int_\Omega\big(R_0(n,p)(n^*p^*-np) 
  + (R_0(n,p)-R_0(n^*,p^*))(1-n^*p^*)\big)\log\frac{np}{n^*p^*}dx \\
  &\le \int_\Omega(R_0(n,p)-R_0(n^*,p^*))(1-n^*p^*)
  \log\frac{np}{n^*p^*}dx.
\end{align*}
We use the local Lipschitz continuity of $R_0$ and $|1-n^*p^*|\le \eta$:
\begin{align*}
  \int_\Omega&(R_0(n,p)-R_0(n^*,p^*))
  (1-n^*p^*)\log\frac{np}{n^*p^*}dx \\
  &\le C\eta\int_\Omega\big(|n-n^*|+|p-p^*|\big)^2 dx
  \le C\eta H^*.
\end{align*}
Collecting these estimates ends the proof of \eqref{5.ineq}.
\end{proof}

\begin{lemma}\label{lem.HL}
Every weak solution to \eqref{1.n}--\eqref{1.bcD} satisfies $H\le C(M^*)(H^*+1)$ and $L\le C(M^*)(1+H^*)$ in $(0,\infty)$, where $L$ is defined in \eqref{5.L}.
\end{lemma}

\begin{proof}
We have $h(n|\bar{n}) = h(n|n^*) + n\log(n^*/\bar{n}) + \bar{n} - n^* \le h(n|n^*) + C(M^*)(1+n)$ and $n\le h(n|n^*) + (e-1)n^*$, so that $h(n|\bar{n}) \le Ch(n|n^*) + C$. Similar estimates hold for $p$ and $D$. Furthermore, $\|\na(V-\bar{V})\|_{L^2(\Omega)}^2\le 2\|\na(V-V^*)\|_{L^2(\Omega)}^2 + 2\|\na(V^*-\bar{V})\|_{L^2(\Omega)}^2 \le 2\|\na(V-V^*)\|_{L^2(\Omega)}^2 + C$. Combining the estimates gives the bound for $H(t)$, and Lemma \ref{lem.cons} shows the bound for $L(t)$.
\end{proof}

We turn to the proof of Theorem \ref{thm.exp2}. Let $H_0:=H^*(0)>0$ (otherwise, the solution coincides with the steady state by uniqueness). We define $L_*:=C(M^*)(H_0+2)$ with $C(M^*)$ from Lemma \ref{lem.HL}; the bounds $m:=m(L_*)$ and $M:=M(L_*)$ from Proposition \ref{prop.unifbd}; the constant $C_1$ from Lemma \ref{lem.rhs} with the values $m$ and $M$; and ${\eta_0}>0$, $\kappa>0$ will be chosen later, depending on $C_1$. All these constants depend only on $H_0$, $M^*$, and the data.

The regularity $n,p,D\in C^0_{\rm loc}([0,\infty);L^2(\Omega))$ and Lemma \ref{lem.Hstar} imply that $t\mapsto H^*(t)$ is absolutely continuous on $[0,\infty)$. We set $\eta_0:=1/(2C_1(C_1+1))$ and $\kappa:=1/(2C_1)$, and we assume that $\eta\le\eta_0$. Note that all constants are fixed before any time interval is considered. Let
\begin{align*}
  S:=\{t\ge 0: H^*(s)\le H_0+1 \mbox{ for all } s\in[0,t]\}.
\end{align*}
Then $0\in S$, and $S$ is relatively closed in $[0,\infty)$, since $H^*$ is continuous. Let $t\in S$. By Lemma \ref{lem.HL}, $L(s)\le C(M^*)(H_0+2)=L_*$ for $s\in[0,t]$, and Proposition \ref{prop.unifbd} gives $m\le n,p,D\le M$ in $\Omega_t$, where $m$ and $M$ do not depend on $t$, because Proposition \ref{prop.unifbd} gives bounds independently of the length of the interval. Therefore, Lemma \ref{lem.rhs} applies on $(0,t)$ with the constant $C_1$. We insert \eqref{5.ineq} into \eqref{5.dHstar} to find for a.e.\ $s\in(0,t)$ that
\begin{align*}
  \frac{dH^*}{ds}\le -\mathcal{D}^* + C_1\eta(\mathcal{D}^*+H^*)
  = -(1-C_1\eta)\mathcal{D}^* + C_1\eta H^*.
\end{align*}
Since $1-C_1\eta\ge 0$, we infer from the inequality $H^*\le C_1\mathcal{D}^*$ from Lemma \ref{lem.rhs} and $\eta\le\eta_0$ that
\begin{align*}
  \frac{dH^*}{ds}\le -\big(C_1^{-1}(1 - C_1\eta) - C_1\eta\big)H^* 
  \le -\kappa H^* \quad\mbox{for a.e. }s\in(0,t).
\end{align*}
Hence, $H^*(s)\le H_0e^{-\kappa s}\le H_0$ for $s\in[0,t]$. In particular, $H^*(t)<H_0+1$, and by continuity, there exists $\tau>0$ such that $H^*<H_0+1$ on $[t,t+\tau]$, i.e., $[0,t+\tau]\subset S$. Thus, $S$ is also relatively open in $[0,\infty)$. Since $[0,\infty)$ is connected, $S=[0,\infty)$, and the above estimates hold for all $t>0$, which means that
\begin{align*}
  H^*(t)\le H_0 e^{-\kappa t}, \quad m\le n,p,D\le M
  \quad\mbox{for all }t>0.
\end{align*}
Finally, because of the density bounds and Taylor expansion, the inequality $\|n-n^*\|_{L^2(\Omega)}^2 + \|p-p^*\|_{L^2(\Omega)}^2 + \|D-D^*\|_{L^2(\Omega)}^2 + \|\na(V-V^*)\|_{L^2(\Omega)}^2\le CH^*(t)$ finishes the proof{, with $\kappa_1:=\kappa/2$}.

\begin{remark}[Three space dimensions]\label{rem.3d}\rm 
In three space dimensions, the Gagliardo--Niren\-berg inequality \eqref{1.GN} holds with $\theta=3/2-3/q$. To estimate the drift term, we apply H\"older's inequality with $1/q+1/r+1/2=1$. Inserting $q=2r/(r-2)$ gives $\theta=3/r$. Therefore, we need the elliptic regularity $V\in W^{1,r}(\Omega)$ with $r>3$. This is possible: Under the geometric conditions of \cite{DKR15} (in particular, $\Gamma_D$ and $\Gamma_N$ intersect in the plane at an angle not larger than $\pi$), the elliptic regularity improves to $r>3$. The norm
\begin{align*}
  \|\na V\|_{L^\infty(0,T;L^r(\Omega))}
  \le C + C\|n-p-D+A\|_{L^\infty(0,T;L^{3r/(3+r)}(\Omega))}
\end{align*}
is bounded if the densities are bounded in $L^\infty(0,T;L^s(\Omega))$ with $s>3/2$. Thus, we need more than the $L\log L$ bounds from the free energy inequality. {\em Assuming} the regularity in $L^s(\Omega)$ with $s>3/2$, the $L^\infty(\Omega_T)$ bounds of Lemma \ref{lem.up} and the uniqueness proof of Theorem \ref{thm.ex} hold in three space dimensions.

The problem is to achieve uniform bounds for the densities better than $L\log L$. In particular, we cannot extend the proof of $L^2(\Omega)$ bounds in Lemma \ref{lem.L2H1} to three dimensions. The proof works in two dimensions, since the $L\log L$ bounds are sufficient to estimate the drift term via the logarithmic-type Gagliardo--Nirenberg inequality for $d=2$. In fact, if the densities are bounded in $L^2(\Omega)$ (uniformly in time), we can apply the proofs of Theorems \ref{thm.exp1}--\ref{thm.exp2} on the exponential time decay. However, the proof of $L^2(\Omega)$ bounds ($L^s(\Omega)$ bounds with $s>3/2$ are sufficient) in three dimensions is currently an open problem. 
\end{remark}

%%%%%%%%%%%%%%%%%%%%%%%%%%%%%%%%%%%%%%%%%%%%%%%%%%%%%%%%%%%%

\end{document}